\documentclass[12pt, draft]{article}
\usepackage{amssymb,amsmath,eucal,latexsym,amsthm}
\usepackage{indentfirst}
\usepackage{textcomp}
\usepackage{enumerate}
\usepackage{verbatim}

\usepackage[dvips]{graphicx}
\theoremstyle{plain}
\newtheorem{theorem}{Theorem}[section]
\newtheorem*{theom}{Main Theorem}

\newtheorem{corollary}{Corollary}[section]
\newtheorem{proposition}{Proposition}[section]
\newtheorem{lemma}{Lemma}[section]

\theoremstyle{definition}
\newtheorem{definition}{Definition} 

\newtheorem*{examples}{Examples}
\theoremstyle{remark}
\newtheorem*{remark}{Remark}
\newtheorem*{remarks}{Remarks}

\DeclareMathOperator{\Hol}{Hol}
\DeclareMathOperator{\rad}{rad}
\newcommand{\R}{\mathbb{R}}
\newcommand{\bC}{\mathbb{C}}
\newcommand{\N}{\mathbb{N}}

\newcommand{\D}{\mathbb{D}}
\newcommand{\e}{\varepsilon}
\DeclareMathOperator{\Zero}{Zero}
\DeclareMathOperator{\Har}{Harm}

\DeclareMathOperator{\ac}{ac}
\DeclareMathOperator{\sic}{sc}
\DeclareMathOperator{\ud}{ud}
\DeclareMathOperator{\conv}{conv}
\DeclareMathOperator{\dist}{dist}

\DeclareMathOperator{\supp}{supp}
\DeclareMathOperator{\reg}{reg}

\newcommand{\mc}{\mathcal}
\renewcommand{\d}{{\rm \, d}}
\renewcommand{\Re}{{\rm Re \,}}
\renewcommand{\phi}{\varphi}

\DeclareMathOperator{\loc}{loc}

\title{Zero (sub)sets \\for spaces of holomorphic functions\\ 
and (sub)harmonic minorants
}
\author{
B.\,N.~Khabibullin\thanks{Supported by the RFSB grant No.~03--01--00033,
and by the Russian programm ``State support of leading scientific schools'', project No.~1528.2003.1.}
\\
\\
{\it {\footnotesize Department of Mathematics, Bashkir State University,}}\\
{\it {\footnotesize Frunze Str., {\rm 32}, Ufa, Bashkortostan, {\rm 450074},Russia;}}\\
{\it {\footnotesize Institute of Mathematics with CS of UCC of RAS,}}\\
{\it {\footnotesize   Chernyshevskii Str., {\rm 112}, Ufa, Bashkortostan, {\rm 450077}, Russia}}\\
{\footnotesize E-mail: khabib-bulat@mail.ru}\\
\footnotesize Web-site: www.bashedu.ru/khabib-bulat
}

\date{}
\begin{document}

\maketitle
\begin{abstract}
Let $\Lambda =\{ \lambda_k \}$ be a sequence of points in a domain
$\Omega $ of the complex plane $\bC$.
We obtain various general conditions in terms of the balayage and Green's functions under which the sequence $\Lambda$
is the zero (sub)set for weighted spaces of holomorphic functions in $\Omega$.
In particular, we consider the classical space $A^{-p}$ with $p>0$,
i.\,e., the set of functions $f$ holomorphic in the unit disk $\D$ satisfying 
$\sup_{z\in \D} |f(z)|(1-|z|)^{p}<+\infty$.
\end{abstract}

\tableofcontents

\section*{Introduction}
\addcontentsline{toc}{section}{Introduction}

Let $\Omega$ be a domain in the complex plane $\bC$ with the boundary $\partial \Omega$. 
Denote  by $\Hol(\Omega)$  the space of all holomorphic functions in $\Omega$. We are concerned with finite or infinite sequences $\Lambda =\{ \lambda_k \}$, $k=1,2, \dots$ of not necessarily distinct points from the domain $\Omega$, without limit points in $\Omega$. 
Let $n_{\Lambda}$ be an integer-valued {\it counting measure\/} of sequence $\Lambda$, defined by 
\begin{equation}\label{df:nLa}
n_{\Lambda}(S):=\sum_{\lambda_k \in S} 1	, \quad S\subset \Omega .
\end{equation}
The function 
\begin{equation}\label{df:dLa}
n_{\Lambda}(z):=n_{\Lambda}\bigl( \{z  \}\bigr),
\quad z \in \Omega,
\end{equation}
is called the {\it divisor} of sequence $\Lambda$. In our paper the sequence $\Lambda$ coincides with a sequence $\Gamma =\{\gamma_n\}$
(or is equal to $\Gamma$, or $\Lambda =\Gamma$)
if and only if $n_{\Lambda}(z)=n_{\Gamma}(z)$ for every $z\in \Omega$; $\Lambda \subset \Gamma$
means $n_{\Lambda}(z)\leq n_{\Gamma}(z)$ for all $z\in \Omega$.

Given $f\colon A\to B$ and $b\in B$ , we write {\it $f\equiv b$ on $A'$} if $f$ is identically equal to $b$ on $A' \subset A$; in the opposite case,  $f\not\equiv b$ on $A'$. For a subset $A' \subset A$,
denote by $f \bigm|_{A'}$ the restriction of $f$ to $A'$.

Let $A,B\subset [- \infty, +\infty ]$.
A function $f\colon A\to B$ is {\it increasing\/} ({\it decreasing} resp.) if, for any $x_1, x_2\in A$, the inequality $x_1\leq x_2$ implies  the inequality $f(x_1)\leq f(x_2)$
\; ( $f(x_1)\geq f(x_2)$ resp.). 
  
``Positive'' (``negative'' resp.) means ``$\geq 0$'' (``$\leq 0$'' resp.).

Let $f\in \Hol(\Omega)$, $f\not\equiv 0$ on $\Omega$. We denote by $\Zero_f $ the sequence of zeros, counting multiplicities, of the function $f$. The sequence $\Lambda$ is {\bf zero set} for a class $H\subset \Hol(\Omega)$
(we write $\Lambda \in \Zero(H)$) if and only if there exists a function $f\in H$ such that $\Lambda = \Zero_f$.
A function $f\in \Hol(\Omega)$ {\it vanish on\/} $\Lambda$ if and only if  $\Lambda \subset \Zero_f$ (we write $f(\Lambda)=0$). 
The sequence $\Lambda$ is {\bf zero subset} for a class $H\subset \Hol(\Omega)$  if and only if there exists a nonzero function $f\in H$ which vanish on $\Lambda$. 

By $\D$ denote the {\it unit disk\/} $\{ z\in \bC \colon |z|<1\}$.

By $m$  denote the {\it Lebesgue measure\/} on $\bC$.

We write $D\Subset \Omega$ if the closure $\overline D$ of $ D\subset \Omega$ is a compact subset in $\Omega$, i.\,e.,
the set $D$ is a relatively compact subset of $\Omega$. 

For domain $D \subset \bC$,  denote by $g_D(\cdot , z)$ 
the {\it extended Green's function\/} for $D$ with a pole at $z\in D$
\cite[5.7.4]{HK}, i.\,e., $g_D(\zeta , z )\equiv 0$ for all $\zeta \in \bC \setminus \overline D$, and $g(\zeta , z)$ is a subharmonic function of $\zeta \in \bC \setminus \{ z \}$. 

For a real-valued function $M$ on an interval of the real axis $\R$,  denote by $M_-'$ the {\it left-hand derivative\/} of $M$.

Given  $M \colon \Omega \to [-\infty , +\infty]$, we put
\begin{equation}\label{sp:HM}
\Hol(\Omega; M){:=}
\left\{ f\in \Hol(\Omega) \colon 
\sup_{z\in \Omega} \frac{|f(z)|}{\exp M(z)}<+\infty
\right\} \, .
\end{equation}

A function $M \colon \D \to \R$ is called {\it radial\/} if
$M(z)=M(|z|)$ for all $z\in \D$.
Given $z\in \bC$ and $t \in \R$, denote by $D(z,t)$ the open disk of radius $t$, centered at $z$. If $t\leq 0$, then $D(z, t)$ is the {\it empty set\/} $\varnothing$. By definition $D(t):=D(0,t)$.

{\it Hereinafter by\/ $\Omega$ denote a domain in\/ $\bC$ containing
the origin, i.\,e.,\/ $0\in \Omega$}.

By $\mc U_0^d(\Omega )$ denote the class of all subdomains $D\Subset \Omega$ such that each domain $D \in \mc U_0^d(\Omega )$ is {\it the union of finite number of disks\/} $D(z,t)\Subset \Omega$, and $D$ contains the origin, i.\,e., $0\in D$. 

Let $\Omega_0 \neq \varnothing$ be a relatively compact subdomain of $\Omega$, and the domain $\Omega_0$ contains the origin.
By $\mc U_0^d(\Omega ; \Omega_0 )$ denote the subclass of 
$\mc U_0^d(\Omega )$ such that each domain $D \in \mc U_0^d(\Omega ; \Omega_0 )$ {\it includes the domain\/} $\Omega_0$; $\mc U_0^d(\Omega ; \varnothing ) :=\mc U_0^d(\Omega )$.    

The following Theorems \ref{th:o} and \ref{th:0} 
illustrate a some of our results.

\begin{theorem}\label{th:o} Let $M \colon \D \to [0, +\infty )$ be a radial positive function that is continuous at\/ $0$. Suppose that 
$M\bigm|_{[0, 1)}$ is the increasing convex function of\/ 
           $\log$ on $(0, 1)$, i.\,e., 
the superposition $M \circ \exp$ is convex function on $(-\infty, 0)$. 
Under the condition 
\begin{equation}\label{cond:inM}
\int_{r}^{1^-}(1-t) \d M(t) =O(1-r) \quad \text{as \; $r\to 1^-$}, 
\quad t,r \geq 0,
\end{equation}

the following three statements 
hold:
\begin{enumerate}
	\item[{\rm ($\mathrm Z_{r}$)}] 	A sequence $\Lambda =\{ \lambda_n\}$, $0\notin \Lambda$, is\/ {\bf zero set} for the space $\Hol (\D ; M)$ 	if and only if 	there are  constants  $a <1$ and $C$ such that 
the inequality  
\begin{equation}\label{in:maz1}
\sum_n g_D(\lambda_n , 0)\leq
\int_0^{1^-} \left(\frac{1}{2\pi} \int_0^{2\pi} g_D (te^{\theta} , 0 )\d \theta \right)\d \bigl( tM_-'(t)\bigr) +C
\end{equation}
holds for each 
domain $D \in \mc U_0^d \bigl(\D ; D(a ) \bigr)$; 

\item[{\rm ($\mathrm S_{r}$)}]	A sequence\/ $\Lambda$ in\/ $\D$ is {\bf zero subset} for the space $\Hol(\D; M)$ if and only if the sequence $\Lambda$ is a {\bf zero set} for the space $\Hol(\D; M)$; 

\item[{\rm ($\mathrm M_{r}$)}] If $f=g/q$ is a meromorphic function on\/ $\D$, and $g, q\in \Hol (\D;M )$, then there are holomorphic functions $g_0, q_0\in \Hol (\D;M )$  {\bf without common zeros} such that $f=g_0/q_0$ on $\D$. 
\end{enumerate}
\end{theorem}

Theorem \ref{th:o} is a very special case of Theorem \ref{th:MB} $\Longleftarrow$ 
Corollary \ref{th:Mc} $\Longleftarrow$ the Main Theorem (see exact wordings
of these results in Subsection \ref{subsec:12}).
Corollary \ref{th:Mc} is established for spaces $H(\Omega; M)$ in case that $\Omega$ is a simply connected domain, and $M \in SH(\Omega)$.  
In Theorem \ref{th:MB} we consider generally speaking a case of non-radial subharmonic function $M$ on $\D$.

In our article a main model class $\Hol (\D ; M)$ is the space    
$A^{-p}$, $p\geq 0$, with 
$$
M(z)\equiv M(|z|)\equiv p \log \frac{1}{1-|z|} \, , \quad z\in \D ,
$$
i.\,e., the set of functions $f\in \Hol (\D )$  satisfying 
\begin{equation*}
|f(z)| \leq C_f \left(\frac1{1-|z|}\right)^{p} \, ,
\quad \forall z\in \D ,	
\end{equation*}
where $C_f$ is a constant dependent on $f$.
For $p =0$ this space is the space $H^{\infty}$ of bounded holomorphic functions in $\D$. The classical Nevanlinna theorem gives
precise geometric information about the zero sets for $H^{\infty}$:
{\it a sequence\/ $\Lambda$ is zero set for\/ $H^{\infty}$ if and only if\/ $\sum_{\lambda_k \in \D}\bigl(1-|\lambda_k|\bigr) <+\infty$} (the classical Blaschke condition), and 
{\it class of zero sets for\/ $H^{\infty}$ coincides with the class of of zero subsets for $H^{\infty}$}.
B.~Korenblum introduced in the work \cite{Korenblum}  a notion
of density, in a certain sense generalizing the classical Blaschke condition, and
found a complete geometric description of zero sets for the algebra $A^{-\infty}=\cup_{0\leq p}  A^{-p}$.
E.~Beller \cite{Beller} proved that   for any $p >0$, the class of zero sets for $A^{-p}$ coinsides with the class of zero subsets for $A^{-p}$.  
D.~Pascuas \cite{Pascuas} and J.~Bruna and X.~Massaneda \cite{BM}
generalized the result of B.~Korenblum to weighted algebras 
$A^{\lambda}(\D):=\bigl\{ f\in \Hol(\D )\colon \log |f(z)| \leq C_f \lambda (|z|)\bigr\}$	
 where $\lambda >0$ is a ``slowly increasing'' function, and $C_f$
 is a constant.
 
In \cite{Seip94}--\cite{Seip} K.~Seip evolved  the method of Korenblum
and obtained a similar complete description of zero set for spaces 
$A_+^{-p}= \cap_{p' >p} A^{-p'}$ and $A_-^{-p}= \cup_{p' <p} A^{-p'}$.
The joint book of H.~Hedenmalm, B.~Korenblum, and K.~Zhu \cite[Ch.~4]{HKZ} contains a detailed analysis of this results together with their improved interpretations. 

D.~Luecking gave in \cite{Luecking} a criterion of zero sets for $A^{-p}$ in terms of existence of harmonic majorant for special test function constructed by $\Lambda$.  
A development of Luecking's method to weighted spaces on $\D$ 
with ``slowly increasing'' weight was recently proposed in \cite{BKN}.

Following to  \cite{BH}, we will refer to the spaces  $A^{-p}$ as {\it uniform 
Bergman spaces}.

In our article we  establish necessary and sufficient conditions of zero sets for weighted classes of holomorphic functions and their corollaries for the uniform Bergman spaces  in terms of  balayage of  measures and functions (see Theorems \ref{th:i} and \ref{th:l} in Section \ref{sec:6}), and in terms of Green's functions.

\begin{theorem}\label{th:0} Let $\Lambda=\{ \lambda_k\}$, $k=1, 2, \dots$, $0\notin \Lambda$, be a sequence of points in $\D$ 
and\/ $0\leq p < +\infty$. The following three statements are equivalent:
\begin{enumerate}[{\rm (i)}]
\item\label{i0} $\Lambda$ is a zero set for $A^{-p}$;
\item\label{ii0} There exist constants $a<1$ and $C$ such that 
for any $D\in \mc U_0^d \bigl(\D ; D(a ) \bigr)$,
\begin{equation}\label{t:g0}
\sum\limits_k         
g_D(\lambda_k , 0) \leq 
p \int_0^1 \left(\frac1{2\pi}\int_0^{2\pi}
g_D(te^{i\theta},0)\d \theta \right)\, \frac{dt}{(1-t)^2} +C \, ; 
\end{equation}
\item\label{ivl0} 
There are constants $a <1$ and $C$ such that for any $D\in \mc U_0^d\bigl(\D ; D(a) \bigr)$, 
\begin{multline}\label{t:L0}
 \sum_k   
\bigl( 1-|\lambda_k|^2\bigr)^2
\left(\frac{1}{\pi}\int_{\D}g_D(\zeta ,0) \, \frac{\d m(\zeta )}{|1-{\lambda_k}\overline{\zeta}|^4}\right) \\
\leq p \int_0^1 \left(\frac1{2\pi}\int_0^{2\pi}
g_D(te^{i\theta},0) \d \theta \right)\, \frac{\d t}{(1-t)^2} +C.
\end{multline}
\end{enumerate}
\end{theorem}

The author thanks Daniel H.~Luecking for sent article \cite{Luecking}, and H\r{a}kan Hedenmalm for sent book 
\cite{HKZ}.

\section{Main notions and results}\label{sec:1}

\setcounter{equation}{0}

Let $\Omega$ be a domain in $\bC$, $0\in \Omega$, and let $S$ be a subset in $\Omega$. We write $C(S)$ for the space of all continuous real-valued functions on $S$. 

By $\mathcal M (S)$ denote the set of all real-valued Borel 
(Radon) measures 
on $S$ (on $C(S)$); by  $\mathcal M^+ (S)$ denote the subcone of $\mathcal M (S)$ consisting of positive Borel measures; by $\mc M^+_{\ac}  (S)$ denote the subset of $\mathcal M^+ (S)$ consisting of measures that are absolutely continuous with respect to Lebesgue measure $m$.

Let $\mu \in \mathcal M (\Omega)$. 
Denote by $\supp \mu$ the {\it support\/} of $\mu$. We say that a measure $\mu \in \mathcal M^+ (\Omega)$ is {\it concentrated in\/} a subset $S\subset \Omega$ if $\mu (\Omega \setminus S)=0$.
Given Borel set $B\subset S$, we denote by $\mu \bigm|_B$
the {\it restriction of\/ $\mu$ to\/ $B$}.
    
Given $\nu \in \mathcal M (\Omega)$, $z\in \Omega$, $t \geq  0$, 
we write
\begin{equation}\label{df:nl}
\nu (z, t):=\nu \bigl( D(z,t)\bigl) 
\text{\;  if $D(z, t) \subset \Omega$}, 
\quad \nu^{\rad} (t):=\nu (0, t). 	
\end{equation}
We hope that last notation will not create confusion  with   notion 
\eqref{df:dLa} for the divisor of sequence $\Lambda$, 
defined by the counting measure \eqref{df:nLa}, and with notation
$\d \nu (\zeta)$ which means that the  variable of integration is 
$\zeta $.    

By $\Har(\Omega)$ denote  the space of all harmonic functions on $\Omega$,
and  by $SH (\Omega)$ denote the cone of all subharmonic functions  on $\Omega$.
The function $\equiv-\infty $ on $\Omega$ belongs to $SH (\Omega)$. Besides, $-SH(\Omega)$ is  the cone of all superharmonic functions  on $\Omega$.
For $u \in SH(\Omega)$, $u\not\equiv -\infty$, we denote by\footnote{Here $\Delta$ is the Laplace operator which acts in sense of distribution theory.}
$\nu_u:=\frac1{2\pi}\Delta u$ the {\it Riesz measure\/} of $u$.

\subsection{Subharmonic kernels}

Let $B$ be a Borel subset of $\Omega$, and $\nu \in \mc M^+(\Omega)$. By definition, $L^1(B,\d \nu)$ is the set of all functions $q\colon B \to [-\infty , +\infty]$ that are integrable with respect to the restriction of $\nu$ to $B$, i.\,e., $\int_B|q| \d \nu <+\infty$.

\begin{definition}\label{df:subk} Let $B$ be a Borel subset of $\Omega$. Let 
\begin{equation*}
h \colon (\zeta , z)
\longrightarrow  \R, \quad (\zeta , z)\in B\times \Omega,
\end{equation*}
be a Borel-measurable function which is locally bounded. Suppose, for every fixed point $\zeta \in B$, the function 
$h (\zeta ,\cdot )$ is {\it harmonic on\/} $\Omega$;
then the function 
\begin{equation}\label{rep:subk}
k (\zeta , z):=
\log |\zeta -z|+h (\zeta , z), \quad (\zeta , z)\in B\times \Omega,	
\end{equation}
is {\it a subharmonic kernel on\/} $B\times \Omega$
({\it supported by $B$ with the harmonic component\/} $h$). 

Suppose a $\nu \in \mc M^+(\Omega)$
is {\it concentrated in\/} $B$, i.\,e., $\nu (\Omega\setminus B)=0$.
A subharmonic kernel $k$ on $B\times \Omega$ is  {\it suitable for\/} $\nu$
if  for any $z\in \Omega$ there are a a subdomain $D_z\ni z$ of $\Omega$  and a function $q \in L^1((\Omega \setminus D_z)\cap B, \d \nu)$ such that 
\begin{equation*}
\bigl| k(\zeta , w )\bigr| \leq  q (\zeta), \quad \forall \zeta \in (\Omega \setminus D_z)\cap B, \quad \forall w \in D_z.
\end{equation*}
\end{definition}

\begin{examples} Consider some frequently occurring subharmonic kernels.
\begin{enumerate}
\item[{\bf 0.}] The function $\log |\zeta -z|$
is a suitable subharmonic kernel on $\Omega\times \Omega$ for all positive Borel mesures with compact support  in $\Omega$.
\item[{\bf 1.}] If the domain $\Omega$ possesses a  Green's function $g_{\Omega}$
($\partial \Omega$ is non-polar), then the functions $-g_\Omega (\zeta , z)$  and $-g_\Omega (\zeta , z)+\log |\zeta|$ are  subharmonic kernels on $\Omega\times \Omega$ and on $(\Omega\setminus \{ 0\})\times \Omega$
respectively. These kernels are suitable for measures $\nu \in \mc M^+(\Omega)$ satisfying
(see \cite[Theorem 4.5.4]{Ransfordb})
$$
\int_\Omega g_\Omega (\zeta , 0) \d \nu (\zeta )<+\infty .
$$

\item[{\bf 2.}] Let $\Omega=\D$. Here we use the Blaschke factor, a variant of this factor, and the pseudohyperbolic  distance
for $\D$: 
\begin{subequations}\label{n:bbr}
\begin{align}
	B_\zeta ( z)&:=\frac{|\zeta|}{\zeta}\,\frac{\zeta -z}{1-\overline{\zeta}z}\, , \; \zeta \in \D \setminus \{ 0 \}, z\in \D ,  
\text{ but } |B_0 ( z)|:=|z|.	\tag{\ref{n:bbr}B}\label{n:bbrB}
	\\
	{\overline B}_{\zeta} ( z)&:=\frac{\overline{\zeta}(\zeta -z)}{1-\overline{\zeta}z}=|\zeta| \, B_\zeta ( z), \quad 
	\zeta \in \D , \; z\in \D , \tag{\ref{n:bbr}\=B}\label{n:bbroB}
	\\
	\rho (\zeta , z)&:=\left| \frac{\zeta -z}{1-\overline{\zeta}z} \right|
=\bigl| B_\zeta ( z)\bigr|=\frac{1}{|\zeta|} \, \bigl|\overline B_\zeta ( z)\bigr|.
\tag{\ref{n:bbr}$\rho$}\label{n:bbrr}
\end{align}
\end{subequations}

The following functions are subharmonic kernels:
\begin{enumerate}
\item[({$\mathrm B_0$})] The subharmonic {\it Blaschke kernel\/}	on $\D \times \D$
is the function $b_1(\zeta , z):=\log |B_{\zeta}(z)|=\log |\rho (\zeta , z)|=-g_{\D}(\zeta , z)$ which is suitable for measures $\nu \in \mc M^+(\D )$ satisfying
$\int_0^{1^-} (1-t)	\d \nu^{\rad}(t)<+\infty$.
Similarly, $\overline b_1(\zeta , z):=\log |\overline B_{\zeta}(z)|=-g_{\D}(z,\zeta)+\log |\zeta|$
is a subharmonic kernel	on $(\D \setminus \{ 0\} )\times \D$
 which   
is suitable for the same measures $\nu \in \mc M^+(\D )$.
	\item[({ ${\overline {\mathrm D}}_p$})] For an integer $p\geq 0$, the subharmonic {\it Dzhrbashyan's kernel of  genus\footnote{For $q>p$, 
	by definition, $\sum\limits_{k=q}^p \dots :=0$, $\prod\limits_{k=q}^p \dots :=1$. Similarly, $\sum\limits_{k\in \varnothing} \dots :=0$, $\prod\limits_{k\in \varnothing} \dots :=1$.}\/} $p$  supported by $\D \setminus \{ 0\}$ is the function
(see \cite{Djrbashian}--\cite{Colwell})
\begin{subequations}
\begin{align*}
\overline d_p(\zeta , z)&:=\log 	\bigl| \overline B_\zeta ( z)\bigr| + \sum_{k=1}^p\frac1{k}\, \Re
(1-\overline B_\zeta ( z))^k \\&=\log \left|\frac{\overline{\zeta}(\zeta -z)}{1-\overline{\zeta}z}\right|+\sum_{k=1}^p\frac1{k}\, \Re 
\left(\frac{1-|\zeta|^2}{1-\overline{\zeta}z}\right)^k 
\end{align*}
\end{subequations}
which coincides with $-g_{\D} (\zeta , z)+\log |\zeta|$ for $p=0$.
This kernel is suitable for measures $\nu \in \mc M^+(\D )$ satisfying
\begin{equation}\label{c:zD}
\int_0^{1^-} (1-t)^{p+1}	\d \nu^{\rad} (t) +
\int_0^{1/2}  \frac{\nu^{\rad} (t)}{t} \d t <+\infty. 
\end{equation}

\item[($\mathrm H_2$)] The subharmonic {\it Horowitz's kernel\/} (see \cite{Horowitz}, \cite{Colwell}) on $\D \times \D$ is the function
$h_2(\zeta , z):=\log 	\bigl| 1-(1- B_\zeta ( z))^2\bigr|$.
This kernel is suitable for measures $\nu \in \mc M^+(\D )$ satisfying
$\int_0^{1^-} (1-t)^{2}	\d \nu^{\rad} (t)  <+\infty$. 

\item[($\mathrm B_s$)] For $0< s \leq 6$, the subharmonic {\it Beller's kernel\/}  on $\D \times \D$ 
is the function   $b_s(\zeta , z):=\log 	\bigl| 1-(1-B_\zeta ( z))^s\bigr| $
(see \cite{Beller}, \cite{Colwell})
which coincides with the Horowitz's kernel for $s=2$ and with the Blaschke kernel for $s=1$. This kernel is suitable for measures $\nu \in \mc M^+(\D )$ satisfying
$\int_0^{1^-} (1-t)^{s}	\d \nu^{\rad} (t) <+\infty$.

\item[(${\overline{ \mathrm B}}_s$)] For $s\geq 1$, the subharmonic {\it Bomash's kernel\/}  supported by $\D \setminus \{ 0\}$ is the function  (see \cite{Bomash}, \cite{B-H})
\begin{equation*}
{\overline b}_s(\zeta , z):= \log 	\bigl| 1-(1-\overline B_\zeta ( z))^s\bigr|.
\end{equation*}
This kernel is suitable for measures $\nu \in \mc M^+(\D )$ satisfying
\eqref{c:zD} for $p=s-1$. 
Below we will use a special case of subharmonic Bomash's kernel with $s=2$:  
\begin{subequations}\label{df:HB}
\begin{align}
\hspace{-2mm}
\overline b_2(\zeta , z):=\log 	\bigl| 1-(1-\overline B_\zeta ( z))^2\bigr|=
\log \left|1-\left(\frac{1-|\zeta|^2}{1-\overline{\zeta}z}\right)^2 \right|& 
\label{df:HBa}
 \\
\! \! \! =\log 	\left(\bigl| \overline B_\zeta ( z)\bigr| \bigl|2- \overline B_\zeta ( z)\bigr|\right)
=\log \frac{|\zeta||\zeta -z|  \bigl|2-|\zeta|^2-\overline{\zeta}z\bigr|}{|1-\overline{\zeta}z |^2}\, .& 
\label{df:HBb}
\end{align}
\end{subequations}
This kernel is suitable for measures $\nu \in \mc M^+(\D )$ satisfying
\begin{equation}\label{c:Bok}
\int_0^{1^-} (1-t)^{2}	\d \nu^{\rad} (t) +
\int_0^{1/2}  \frac{\nu^{\rad} (t)}{t} \d t <+\infty. 
\end{equation}

\item[($\mathrm K_1$)] The subharmonic {\it Korenblum's kernel\/} on $(\D \setminus \{ 0\})\times \D$ (see  \cite{Korenblum}, \cite{Colwell}) is the function  
$k_1(\zeta , z):=\log 	\bigl| B_\zeta ( z)\bigr| + 
\log \dfrac{1}{|\zeta|}\, \Re 
\dfrac{\zeta /|\zeta| +z}{\zeta /|\zeta| -z} \,.  
$
This kernel is suitable for measures $\nu \in \mc M^+(\D )$ satisfying
\eqref{c:Bok}.
\end{enumerate}

Note that for each measure $\nu \in \mc M ( \D )$
we can select  subharmonic kernel (see Dzhrbashyan's factorization theory in \cite{Djrbashians}, \cite{Djrbashianss}, \cite{Colwell}) which is  suitable for $\nu$ and optimal in a certain sense.
 
\item[{\bf 3.}] Let $\Omega=\bC$, $r_0>0$. The following functions are subharmonic kernels:   
\begin{enumerate}
	\item[({ $\mathrm E_q$})] For an integer $q\geq 0$, the subharmonic {\it Hadamard--Weierstrass kernel of genus\/ $q$} on $(\bC \setminus \{ 0\}) \times \bC$	is the function
$$
e_q(\zeta , z):=\log \Bigl|1-\frac{z}{\zeta}\Bigr|+\sum_{k=0}^q \frac1{k}\, \Re \frac{z}{\zeta} \, .
$$
This kernel is suitable for measures $\nu \in \mc M^+(\bC )$ satisfying
$$
\int_0^1 \frac{\nu^{\rad} (t)}{t^{q+1}} \d t
+\int_1^{+\infty} \frac{\nu^{\rad} (t)}{t^{q+2}} \d t <+\infty .
$$
\item[({ $\mathrm W$})] Let $\{q_n  \}$  be a sequence of nonnegative 
integer,  and let $\{ r_n \}$ be an increasing sequence of positive numbers, $n\in \N$, $r_1>r_0$. The subharmonic {\it Weierstrass kernel\/} (with respect to these two sequences) on $\bC \times \bC$ is the function $w(\zeta , z):=e_{q_n}(\zeta , z)$ when $r_{n-1}\leq |\zeta|< r_{n}$,
$n\in \N$, and $w(\zeta , z):=\log |\zeta -z|$ when $|\zeta|< r_0$.
This kernel is suitable for measures $\nu \in \mc M^+(\bC )$ satisfying
(see \cite[Theorem 4.1]{HK})
$$
\sum_{n=1}^{\infty}\int_{r_{n-1}}^{r_n} \frac{\nu^{\rad} (t)}{t^{q_n+2}} \d t <+\infty .
$$
\end{enumerate}	
Note also that, for each measure $\nu \in \mc M^+ ( \bC )$,
we can select  subharmonic kernel (see Dzhrbashyan's factorization theory in \cite{Djrbashianse} and the survey \cite{GLO}) which is  suitable for $\nu$
and optimal in a certain sense.
 
\end{enumerate}
\end{examples}

The following Proposition \ref{pr:Riesz} is a global version of Riesz Decomposition Theorem for subharmonic functions \cite{HK}, \cite{Ransfordb}. 
\begin{proposition}\label{pr:Riesz}
Let $k$ be a subharmonic kernel on $B\times \Omega$. 
Suppose a measure $\nu \in \mc M^+(\Omega)$ is concentrated in $B$.
If the kernel $k$ is suitable for $\nu$, then  
the function 
\begin{equation*}
U^{\nu}_k(z):=	\int_\Omega k (\zeta , z) \d \nu( \zeta )
\end{equation*}
is subharmonic on $\Omega$ with the Riesz measure $\nu$.
In particular, every function $M\in SH(\Omega)$ with the Riesz measure $\nu_M=\nu$
can be decomposed as 
\begin{equation}\label{repr:Riesz}
M=	U_k^{\nu}+H, \quad \text{where } H\in \Har (\Omega).
\end{equation}
\end{proposition}

The proof is omitted here. It can be obtained from \cite[Theorem 2.6.5]{Klimek}. 

\subsection{Main results}\label{subsec:12}

We write  $\, L^1_{\loc}(\Omega)\,$  for the set  of 
all  functions $F \colon \Omega \to [-\infty , +\infty]$ that are locally integrable with respect to $m$. 
A sequence $\{  w_n \}$ from $ L^1_{\loc}(\Omega )\,$
is convergent in $L^1_{\loc}(\Omega )\,$  if there exists a 
function $w\in L^1_{\loc}(\Omega )\,$ such that $\lim_{n\to \infty} \int_K |w_n-w|\d m \to 0$ for any compact subset $K\Subset \Omega$.

We denote by $\dist (z, \partial \Omega)$ the Euclidean distance from $z\in \Omega$ to $\partial \Omega$.

In what follows the measure $m^{(r)} \in \mc M^+(\bC )$ is obtained from Lebesgue measure $m$ by restricting it to the disk $D(r)$ and normalizing so that $m^{(r)} \bigl(D(r)\bigr)=1$, i.\,e.,
\begin{equation}\label{df:mr}
m^{(r)}:=	\frac1{m^{\rad}(r)} \, m\bigm|_{D(r)}.
\end{equation}
   
We have
\begin{equation}\label{df:convo}
v(z)\leq \bigl(v*m^{(r)}\bigr)(z), \quad  v\in SH(\Omega),
\end{equation}
if $D(r)\Subset \Omega$ (here and below $*$ stands for the convolution operation).

 Let  $\sigma \colon \Omega \to (0,+\infty)$  
 be a function satisfying
\begin{equation}\label{c:rD}
0<\sigma (z)<\dist (z, \partial \Omega ), \quad \forall z\in \Omega.	
\end{equation}
For $F\in  L^1_{\loc}(D)$, we put
\begin{equation}\label{df:si}
F^{(\sigma)}(z):=\int_{D(\sigma(z))} F(z+w) \d m^{(\sigma(z))}(w), \quad  z\in \Omega.	
\end{equation}

We say that a function $F : \Omega \to [-\infty , +\infty]$
{\it admits a  harmonic minorant\/} ({\it a  subharmonic minorant\/} resp.) 
on $\Omega$ if and only if there is a $h\in \Har (\Omega)$ ($h\in SH (\Omega)$, $h\not\equiv -\infty $ on $\Omega$, resp.) such that $h(z) \leq F(z)$ for all $z\in \Omega$.

For $a\in [-\infty , +\infty]$
(for a function $f \colon X \to [-\infty , +\infty ]$ resp.), as usual, we put 
$a^+:= \max \{ a , 0\}$ ($f^+(x):=\max \{ f (x) , 0\}$ resp.). 

The following theorem is the main general result of our article.

\begin{theom}\label{Th:m} 
Let $M \not\equiv -\infty$ be a subharmonic function with the Riesz measure $\nu_M$ on a domain $\Omega \subset \bC$, $0\in \Omega$.
Suppose also that 
\begin{enumerate}[{\rm i)}]
\item\label{rmi} the function $M$ is bounded on a subdomain $\Omega_1 \Subset \Omega$, and $0\in \Omega_1$;
\item\label{rmii} the measure $\nu_M$ is concentrated in a Borel subset $B\subset \Omega$, and subharmonic kernel $k$ on $B\times \Omega$ is suitable for $\nu_M$;
\item\label{rmiii} for the function 
\begin{equation}\label{df:Qkn}
Q^{\nu_M}_k(z):=\int_B	\bigl(k(\zeta , 0)-k(\zeta , z )\bigr)^+  \d {\nu}_M (\zeta ), \quad z\in \Omega,
\end{equation}
there exists a majorizing function $Q \in L^1_{\loc} (\Omega)$ on\/ $\Omega$, i.\,e., $Q^{\nu_M}_k(z)\leq Q(z)$ at almost all (with respect to Lebesgue measure $m$) points $z\in \Omega$.
\end{enumerate}

Let $u$ be a subharmonic function with the Riesz measure $\nu_u$ on $\Omega$, and $u(0)\neq -\infty$.
Then the following two assertions hold:
\begin{enumerate}
	\item[{\rm (Z)}] Let $\Omega_0$
	be a relative compact subomain of\/ $\Omega_1$ containing the origin.
	If there exists  a constant  $C$ such that the inequality
\begin{equation}\label{in:ma}
\int g_D(\zeta , 0) \d \nu_u(\zeta)\leq
\int g_D (\zeta , 0 )\d \nu_M (\zeta )+C 
\end{equation}
holds for each domain $D \in \mc U_0^d ( \Omega ; \Omega_0 )$,
then, for any $\sigma \in C(\Omega)$ satisfying\/ \eqref{c:rD},
there exists a\/ {\bf harmonic function} $h$ on $\Omega$ such that
\begin{equation}\label{e:estM}
u(z)+h(z)\leq  M^{(\sigma)}(z)+Q^{(\sigma)}(z), \quad \forall z\in \Omega. 	
\end{equation}
\item[{\rm (S)}] If the difference $M-u$ admits\footnote{If $M(z)=u(z)=-\infty$, then we put $(M-u)(z)=+\infty$.} a\/ {\bf subharmonic minorant} on\/ $\Omega$, then, for any $\sigma \in C(\Omega)$ 
satisfying\/ \eqref{c:rD}, the function 
\begin{equation}\label{rep:MQu}
M^{(\sigma)}+Q^{(\sigma)}-u	
\end{equation}
 admits a\/ {\bf harmonic minorant} on\/ $\Omega$. 
\end{enumerate}
\end{theom}

An important consequence from the Main Theorem is

\begin{corollary}\label{th:Mc}
Let $\Omega \subset \bC$ be a simply connected domain containing the origin. Let $M \not\equiv -\infty$ be a subharmonic continuous function on $\Omega$ with the Riesz measure $\nu_M $.  
Suppose that the conditions\/ {\rm \ref{rmi})--\ref{rmiii})}
of the Main Theorem are fulfilled.
Then the following three assertions hold:
\begin{enumerate}
	\item[{\rm ($\mathrm Z_{\sic}$)}] 	Let $\Lambda =\{ \lambda_n\}$ be a sequence in\/ $\Omega$, $0\notin \Lambda$, and let\/ $\Omega_0$ be a relative compact subdomain of\/ $\Omega$ containing the origin.	If there exists  a constant\/  $C$  such that 
the inequality 
\begin{equation}\label{in:maz}
\sum_n g_D(\lambda_n , 0) \leq
\int g_D (\zeta , 0 )\d \nu_M (\zeta )+C  
\end{equation}
holds 
for each domain $D \in \mc U_0^d(\Omega ; \Omega_0 )$,
then, for any\/ $\sigma \in C(\Omega)$ satisfying\/ \eqref{c:rD}, the sequence\/ $\Lambda$ is a {\bf zero set} for\/ $\Hol (\Omega; 
M^{(\sigma)}+Q^{(\sigma)})$.
	\item[{\rm ($\mathrm S_{\sic}$)}]	If a sequence $\Lambda$ in $\Omega$ is {\bf zero subset} for the space $\Hol(\Omega;M)$, 
then, for any\/ $\sigma \in C(\Omega)$ satisfying\/ \eqref{c:rD},
the sequence $\Lambda$ is a {\bf zero set} for the space $\Hol (\Omega; M^{(\sigma)}+Q^{(\sigma)})$. 
\item[{\rm ($\mathrm M_{\sic}$)}] Suppose that\/  
$f=g/q$ is a meromorphic function on\/  $\Omega$ where $g$ and $q$ belong to\/ $\Hol (\Omega;M )$; then, for any\/ $\sigma \in C(\Omega)$ satisfying\/ \eqref{c:rD}, there are holomorphic functions $g_0, q_0 \in \Hol (\Omega;M^{(\sigma)}+Q^{(\sigma)})$ {\bf without common zeros} such that $f=g_0/q_0$ on\/ $\Omega$. 
\end{enumerate}
\end{corollary}

Given $z\in \D$ and $\alpha >0$, the set 
\begin{equation}\label{df:box}
  {\Box}_\alpha(z)
 :=\{ \zeta \in \D \colon |z|-\alpha (1-|z|)\leq |\zeta |, \; 
|\arg \zeta -\arg z|\leq \alpha (1-|z|)\}.
\end{equation}
is called the {\it Carleson box\/} of relative size $\alpha$ with center at $z$.

Given $0<\e <1$ and $M \colon \D \to [-\infty , +\infty)$ , we put
\begin{equation}\label{df:AMe}
A_{M, \e}(z):=\frac1{2\pi}\int_0^{2\pi } M \bigl(z+ \e (1-|z|)e^{i\theta}\bigr) \d \theta	
\end{equation}
if this integral there exists.

The following Theorem \ref{th:MB} is a special version of previous
Corollary \ref{th:Mc} in  case  that we use the subharmonic Bomash's kernel
\eqref{df:HB} from the item (${\overline{ \mathrm B}}_s$) of Examples
in the role of subharmonic kernel $k$ on $(\D \setminus \{ 0\})\times \D$. 
\begin{theorem}\label{th:MB} Let $M\not\equiv -\infty$ be a subharmonic function $M$ on $\D$ with the Riesz measure $\nu_M$. Suppose 
that the function $M$ is  bounded below on every $D \Subset \D$.
Suppose also that the condition 
\begin{equation}\label{cond:mea}
	\int_{0}^{1^-} (1-t)^2 \d \nu_M^{\rad} (t)<\infty 
\end{equation}
holds. By definition, put (see \eqref{df:box})
\begin{equation}\label{df:bMe}
b_M^{[\alpha]}(z):= \frac{1}{(1-|z|)^2} \, \int_{\Box_\alpha(z)}(1-|\zeta|)^2 \d \nu_M (\zeta ).
\end{equation}

Under these conditions and notations, the following three assertions hold:
\begin{enumerate}
	\item[{\rm ($\mathrm Z_{\D}$)}] 	Let\/ $\Lambda =\{ \lambda_n\}$ be a sequence in\/ $\D$, and\/ $a<1$. 	If there exists a  constant\/ $C$ such that the inequality\/ \eqref{in:maz} 
holds for each domain\/ $D \Subset \mc U^{d}_0(\D ; D(a)) $,  
then, for any\/ $0<\e <1$, the sequence\/ $\Lambda$
is a {\bf zero set} for the space 
\begin{equation}\label{sp:defq}
\Hol \bigl(\D ; 2A_{M, \e}-M+C_{\e} \cdot b_M^{[\alpha]}\bigr),	
\end{equation}
where\/ $C_{\e}$ is a positive constant dependent only on\/ $\e$, and\/ $\alpha$ is an absolute constant.

	\item[{\rm ($\mathrm S_{\D}$)}]	If a sequence $\Lambda$ in $\D$ is {\bf zero subset} for the space $\Hol(\D; M)$, 
then, for any $0<\e <1$, the sequence $\Lambda$ is a {\bf zero set} for the space \eqref{sp:defq}. 
\item[{\rm ($\mathrm M_{\D}$)}] If $f=g/q$ is a meromorphic function on $\D$, and $g, q\in \Hol (\D;M )$, then, for any $0<\e<1$, there are holomorphic functions $g_0, q_0$  {\bf without common zeros} on $\D$ such that 
$f=g_0/q_0$ on $\D$, and the functions $g_0, q_0$ belong to the space \eqref{sp:defq}. 
\end{enumerate}
\end{theorem}

\section{Bases of general approach}\label{sec:2}

\setcounter{equation}{0}

\subsection{The balayage:\\ 
Arens--Singer and Jensen measures and functions}
Let $H$ be a convex subcone in $SH(\Omega)$, $z\in \Omega$. Let $\delta , \mu
\in \mathcal M^+ (\Omega)$ are measures with compact support on $\Omega$.
We write $\delta \prec_H \mu $ and say that $\mu$ is a {\it balayage of\/ $\delta$
 with respect to\/} $H$ if
\begin{equation}\label{df:bal}
\int h \d \delta \leq \int h  \d \mu , \quad \forall h\in H.
\end{equation}
 
Given $z\in \bC$ , we denote by $\delta_z$ the Dirac measure
at the point $z\in \bC$, i.\,e. $\supp \delta_z =\{ z \}$, $\delta_z (\bC)=1$. 

\begin{definition}[{\rm \cite{KhIs}--\cite{RansfordT}}]\label{df:asj}
A measure $\mu$ with compact support in $\Omega \ni z$  is called
the {\it Jensen\/} ({\it Arens--Singer}, or {\it representing\/} resp.)
{\it measure} for $z$ on $\Omega$ if $\delta_z \prec_H \mu$
where $H=SH(\Omega)$ ($H=\Har (\Omega)$ resp.), i.\,e., if and only if
\begin{equation}\label{df:bal1}
h(z) \leq  \int h  \d \mu,	\; \forall h\in SH(\Omega) \; 
\Bigl(h(z) =\int h  \d \mu, \; \forall h\in \Har (\Omega) \text{\;resp.}\Bigr).
\end{equation}
Denote by $\mc {\mc J}_z(\Omega)$ ($\mc A \mc S_z (\Omega)$ resp.) the class of all
Jensen (Arens--Singer resp.) measures for $z$. If $z=0 \in \Omega$, then
 ${\mc J}(\Omega):={\mc J}_0(\Omega)$, ${\mc A \mc S}(\Omega):={\mc A \mc S}_0(\Omega)$.  
\end{definition}
Evidently, ${\mc J}_z(\Omega) \subset {\mc A \mc S}_z(\Omega)$, and every $\mu \in {\mc A \mc S}_z(\Omega)$ is 
a probability measure, i.\,e.,
\begin{equation}\label{c:prob}
\mu (\Omega)=1.	
\end{equation}
 A subclass of ${\mc J}_z(\Omega)$ is the class $\mc H_z(\Omega)$ ($\mc H_z^{\reg}(\Omega)$ resp.)
of all harmonic measures $\omega_D(z, \cdot )$ for domains
(for regular for the Dirichlet problem domains resp.) 
 $D\Subset \Omega$ at the point $z$. We set $\mc H(\Omega) :=\mc H_0(\Omega)$, $\mc H^{\reg}(\Omega) :=\mc H_0^{\reg}(\Omega)$.   
The {\it potential of a measure\/ $\mu \in {\mc A \mc S}(\Omega)$} is defined to be function
\begin{equation}\label{df:poten}
V_{\mu}(\zeta ):=\int \log |z-\zeta | \d \mu (z)-\log |\zeta | 
=\int \log \Bigl|1-\frac{z}{\zeta}\Bigr| \d \mu (z) ,
\; \zeta \in \bC \setminus \{ 0\} \, .
\end{equation}

Let $\mu \in {\mc A \mc S}(\Omega)$. For any function $u_{\nu} \in SH(\Omega)$ with the Riesz measure $\nu_u$  {\it the generalized Poisson--Jensen formula\/} 
\cite[Proposition 1.2]{Kh03} holds:
\begin{equation}\label{for:PJ}
u_{\nu}(0)=\int_\Omega u_{\nu} \d \mu -\int_\Omega V_{\mu} \d \nu \, .	
\end{equation}
 
\begin{proposition}\label{prop:Vmu} Let $V_{\mu}$ be the potential \eqref{df:poten} of a measure $\mu \in {\mc A \mc S}(\Omega)$, 
and let $k$ be a subharmonic kernel on $B\times \Omega$. 
Then, for every $\zeta \in B \setminus \{ 0\}$, the potential $V_{\mu}$ can be represented in the form
\begin{equation}\label{repr:Vmu}
V_{\mu}(\zeta)=	\int_{\Omega}  k(\zeta , z)  \d \mu (z)-k(\zeta , 0)=
\int_{\Omega} \bigl( k(\zeta , z)-k(\zeta , 0)\bigr)  \d \mu (z).
\end{equation}
\end{proposition}
\begin{proof} 
By \eqref{rep:subk}, for every $\zeta \in B\setminus \{ 0\}$, we have  representations 
\begin{equation}\label{rep:Kz0}
k(\zeta , z) = \log|\zeta -z|+h(\zeta , z), \quad
k(\zeta , 0)=\log|\zeta|+h(\zeta ,0),
\end{equation}
 where $h(\zeta , z)$ is harmonic function of $z\in \Omega$. Hence, for $\mu \in {\mc A \mc S}(\Omega)$, 
 we obtain
\begin{multline*}
\int_\Omega \bigl(k(\zeta , z )-k (\zeta ,0 )\bigr) \d \mu (z)  \\
\overset{\eqref{rep:Kz0}}{=}\int_\Omega \bigl(\log|\zeta -z|+h(\zeta ,z) -\log|\zeta|-h(\zeta , 0)\bigr) \d \mu (z)\\
\overset{\eqref{c:prob}}{=}\int \bigl(\log|\zeta -z| -\log|\zeta|\bigr) \d \mu (z)
+\int h( \zeta , z) \d \mu (z)-h(\zeta ,0)\\
\overset{\eqref{df:poten}, \eqref{df:bal1}}{=}V_{\mu}(\zeta )+ h(\zeta , 0)-h( \zeta , 0)=V_{\mu}(\zeta )
\end{multline*}  
for each $\zeta \in B\setminus \{ 0\}$, as desired.
\end{proof}

\begin{definition}[{\rm \cite{KhIs}--\cite{Sundberg}}]\label{df:asjf} 
A function $V \in SH \bigl( \Omega \setminus \{ 0\}\bigr)$ 
will be called a {\it Arens--Singer} (or {\it representing}) {\it function\/}
on $\Omega$ if this function satisfies the following two conditions:
\begin{enumerate}[{\rm (1)}]
\item There is  a set $K\Subset \Omega$ such that
$V(\zeta )\equiv 0$ for all $\zeta \in \Omega\setminus K$;
\item\label{V2} $\limsup\limits_{\zeta \to 0} \dfrac {V(\zeta )}{-\log |\zeta |} \leq 1$. 
\end{enumerate}
If $V$  satisfies also the condition
\begin{enumerate}
\item[{\rm (3)}]   $V(\zeta )\geq 0$ for all $\zeta \in \Omega \setminus \{ 0\}$,
\end{enumerate}
then we call $V$ a  {\it Jensen function\/} on $\Omega$.
Denote by  ${\mc P_{AS}}(\Omega)$ and ${\mc P_J}(\Omega)$ the classes of all Arens--Singer and Jensen 
functions respectively.
\end{definition}

Evidently, ${\mc P_J}(\Omega) \subset {\mc P_{AS}}(\Omega)$.
The {\it class\/} $\mathcal P_{\mc H}(\Omega)$ \; ( $\mathcal P_{\mc H}^{\reg}(\Omega)$ resp.) of all {\it extended Green functions\/} $g_D(\cdot , 0 )$ for domains $D\Subset \Omega$ 
(for regular for the Dirichlet problem domains $D\Subset \Omega$ 
 resp.)  is a subclass of ${\mc P_J}(\Omega)$.   

\begin{proposition}[{\rm \cite[Duality Theorem]{Kh03}}]\label{pr:mapP} The map\/ $\mc P \colon  \mu \to V_{\mu}$ is an affine\footnote{It means that $\mc P \bigl(\alpha \mu_1+(1-\alpha)\mu_2 \bigr)=\alpha \mc P(\mu_1)+(1-\alpha )\mc P (\mu_2)$ for $0\leq \alpha \leq 1$.} bijection from ${\mc A \mc S}(\Omega)$ onto ${\mc P_{AS}}(\Omega)$,
and from ${\mc J}(\Omega)$ onto ${\mc P_J}(\Omega)$. 
\end{proposition}

Note that, according to Proposition \ref{pr:mapP}, for $V\in {\mc P_{AS}}(\Omega)$ the condition \eqref{V2} is equivalent to the condition (see \cite[Theorem 3.1.2]{Ransfordb})
\begin{enumerate}
\item[{$(2')$}] ${V(\zeta )} \leq {-\log |\zeta |} +O(1)$ as $\zeta \to 0$. 
\end{enumerate}

Let $\Omega_0 \Subset \Omega$ be a domain containing the origin.

Denote by $\mc A \mc S^{\Omega_0}(\Omega)$  (\; $\mc J^{\Omega_0}(\Omega)$ resp.) the set of all Arens--Singer (Jensen resp.) measures $\mu$ on $\Omega$ such that $\Omega_0 \bigcap \supp \mu = \varnothing$. 

Denote by ${\mc P}_{AS}^{\Omega_0}(\Omega)$ (\; ${\mc P}_{J}^{\Omega_0}(\Omega)$ resp.) 
the set of all Arens--Singer (Jensen resp.) functions $V$ that is harmonic on $\Omega_0 \setminus \{ 0\}$
such that (cf. $(2')$)
\begin{equation}\label{li:V}
V(\zeta)=-\log |\zeta|+O(1), \quad \zeta \to 0. 
\end{equation}

\begin{proposition}\label{pr:isao}
The map\/ $\mc P \colon  \mu \to V_{\mu}$ of Proposition\/ {\rm \ref{pr:mapP}} is an affine bijection from ${\mc A \mc S}^{\Omega_0}(\Omega)$ onto $\mc P_{AS}^{\Omega_0}(\Omega)$,
and from ${\mc J}^{\Omega_0}(\Omega)$ onto $\mc P_J^{\Omega_0}(\Omega)$.  
\end{proposition}
The proof of Proposition \ref{pr:isao} follows from Proposition \ref{pr:mapP} 
and \cite[Proposition 1.4, (4)]{Kh03}.

We recall the notion of classical balayage \cite{Landkof}.
Given a bounded regular for the Dirichlet problem domain $\Omega_1 \Subset \Omega$ and a measure $\mu \in \mc M^+( \Omega )$, the {\it classical balayage of\/} $\mu$ {\it relative to\/} $\Omega_1$ is the measure $\mu^{\Omega_1} \in \mc M^+( \Omega)$ defined by  
\begin{equation}\label{df:baly}
	\mu^{\Omega_1} (E):=\mu ( E\setminus \Omega_1)+\int_{\Omega_1} \omega_{\Omega_1}(z, E) \d \mu (z),
	\quad \text{(Borel $E\Subset \Omega$)}.
\end{equation}
In particular, $\Omega_1 \bigcap \supp \mu^{\Omega_1} =\varnothing$,
and $\mu \prec_{SH(\Omega)} \mu^{\Omega_1}$ \cite[Lemma 6.4]{CR01}. 
Hence 
\begin{equation}\label{eq:AMSJ}
\mc A \mc S^{\Omega_1}(\Omega)=\{ \mu^{\Omega_1} \colon 
\mu \in \mc A \mc S(\Omega) \},  \quad 
\mc J^{\Omega_1}(\Omega)=\{ \mu^{\Omega_1} \colon 
\mu \in \mc J(\Omega) \} 	
\end{equation}
for regular subdomains $\Omega_1 \Subset \Omega$ containing the origin.

Let $D \Subset \bC$ be a regular domain containing
the origin. For short, denote by  $\omega_D$ the harmonic measure for 
$D$ at the origin. 

As before, $\delta_0$ is the Dirac measure at $0$.
\begin{proposition}[{\rm \cite[Ch.~IV, \S~1, and Theorem 4.12]{Landkof}}]\label{pr:harmO} Let $ D \Subset \bC$,  $\Omega'\Subset \bC$, and $D\cup \Omega'$ are regular domains containing the origin.
Then 
\begin{equation}\label{eq:harme}
\delta_0^{\Omega'\cup D}=
\omega_{\Omega'\cup D}=(\omega_D)^{\Omega'\cup D}=	
\omega_{\Omega'}-\left(\omega_{\Omega'}\bigm|_{D}\right) +\left(\omega_{\Omega'} \bigm|_{D}\right)^{\Omega'\cup D}.
\end{equation}
\end{proposition}

\subsection{On the existence of (sub)harmonic minorants}\label{ss:exsh}

\begin{proposition}[{\rm cf. \cite[Proposition 7.1]{Kh012}}]\label{pr:dirmaj} 
Let $F\colon \Omega \to [-\infty, +\infty]$ be  a Borel-measurable 
function. If $F$ admits a harmonic minorant, then 
\begin{equation}\label{c:iha}
-\infty < \inf \left\{ \int_\Omega F \d \mu \colon \mu \in {\mc A \mc S}(\Omega) \right\}. 	
\end{equation}
If $F$ admits a  subharmonic minorant, and $F$ is bounded below in a certain neighbourhood of origin, then 
\begin{equation}\label{c:isha}
-\infty < \inf \left\{ \int_\Omega F \d \mu \colon \mu \in {\mc J}(\Omega) \right\}. 	
\end{equation}  
Evidently, the class ${\mc A \mc S}(\Omega)$ in \eqref{c:iha} can be replaced by  ${\mc J}(\Omega)$,  $\mc H(\Omega)$, $\mc H^{\reg} (\Omega)$, and we can substitute 
$\mc H(\Omega)$ or $\mc H^{\reg}(\Omega)$ for ${\mc J}(\Omega)$ in \eqref{c:isha}. 
\end{proposition}
\begin{proof} Let $h$ be a Borel-measurable function on $\Omega$
and $F\geq h$ on $\Omega$. Then  
$\int F  d\mu \geq \int	h \d \mu  $
for every $\mu \in \mc M^+ (\Omega)$. If the minorant $h$ is harmonic, then  
$\int	h \d \mu \geq h(0)$ for all $\mu \in {\mc A \mc S}(\Omega)$ according to \eqref{df:bal1}, and  \eqref{c:iha} is proved. 

Let us assume that $F $ is bounded below on $D(2\e )\Subset \Omega$, and 
the minorant $h$ is subharmonic. There exists a subharmonic function $h_{\e}$ which  coincides with $h$ on $\Omega\setminus D(\e )$ and harmonic on $D(\e)$ \cite[Theorem 2.18]{HK}. Hence  the difference $h_{\e}-C$ is a subharmonic minorant for $F$ if a constant $C$ is sufficiently large.
Therefore, 
$
\int F \d \mu \geq \int	(h_{\e}-C ) \d \mu \geq h_{\e}(0)-C
$ 
for all $\mu \in {\mc J}(\Omega)$ according to
Definition \ref{df:asj}, and  \eqref{c:isha} is proved. 
\end{proof}

We need another treatments of Proposition \ref{pr:dirmaj} for the cases
special functions $F$  representing by the difference 
 \begin{equation}\label{rep:FuM}
F=M-u, \quad u\in SH(\Omega ), \quad u(0)\neq -\infty .	
\end{equation}
Here the functions $M$ and $u$ can  take on the value $-\infty$ at a point $z \in \Omega$ simultaneously
if $M \in SH(\Omega)$ or $M\in L^1_{\loc}( \Omega)$. In this case we will be consider 
that $F(z)=+\infty$.
The following Proposition \ref{pr:dirmajs} is an easy corollary 
of Propositions \ref{pr:dirmaj}, \ref{pr:mapP}, and of generalized Poisson--Jensen formula \eqref{for:PJ}. 

\begin{proposition}
\label{pr:dirmajs} Let $M\colon \Omega \to [-\infty, +\infty]$ be a Borel-measurable function, and $u$ be a subharmonic function on $\Omega$ with the Riesz measure $\nu_u$.

If the function $F :=M-u$ from \eqref{rep:FuM} admits a harmonic minorant 
on $\Omega$, then 
there is a constant $C$ such that
\begin{equation}\label{c:ihah}
\int_{\Omega } V_{\mu} \d \nu_u \leq 
\int_\Omega M \d \mu +C 	
\end{equation}
for all $\mu \in {\mc A \mc S}(\Omega)$ where $V_{\mu}$ is the potential
of $\mu$ (see definition \eqref{df:poten}). Besides, if $M \in SH(\Omega)$, then there is a constant $C$ such that
\begin{equation}\label{c:ihahM}
\int_{\Omega } V \d \nu_u \leq 
\int_\Omega  V \d \nu_M +C 	
\end{equation}
for all Arens--Singer functions $V $ on $\Omega$.

If $F=M-u$ from \eqref{rep:FuM} admits\/ a  subharmonic minorant, and $M$ is bounded below in a certain neighbourhood of origin, then there is a constant $C$ such that \eqref{c:ihah} holds for all $\mu \in {\mc J}(\Omega)$.
Besides, if $M \in SH(\Omega)$, then there is a constant $C$ such that 
\eqref{c:ihahM} holds for all $V \in \mc P_{J}(\Omega)$.
\end{proposition}

Let $\{ u_n \}$ be a sequence of functions subharmonic in $\Omega$ and uniformly locally bounded above (that is, uniformly bounded above on every $G\Subset \Omega$) in $\Omega$. 
We denote by ${\limsup}^{*} u_n$ the upper semicontinuous regularization of the (pointwise) upper limit of $\{ u_n \}$, which is, in turn, a subharmonic function.

The following Theorem is in a sense converse to Proposition \ref{pr:dirmaj}. 
\begin{theorem}[{\rm \cite[Theorem 7.1]{Kh012}}]\label{th:5} Let $H\subset SH(\Omega)$ be a convex cone. Assume that the cone $H$ contains a negative $( \, \leq 0\, )$ function, and let $F \colon \Omega \to [-\infty, +\infty]$ be a function that belongs to $L_{\loc}^1(\Omega)$. 

Let us suppose that for any $K\Subset \Omega$ and any constant $C$
there is a $h\in H$ such that $h\leq C$ on $K$, and at least one of the following conditions holds:
\begin{enumerate}
	\item[$(\mathrm L^*)$] if a sequence $\{ h_n \}$, $h_n\in H$, is locally bounded above, then the function ${\limsup}^*h_n$ belongs to $H$,
\item[{\rm (CL)}] the cone $H$ is sequentially closed in $L^1_{\loc}(\Omega)$.
\end{enumerate}
If 
\begin{equation}\label{i:absb}
-\infty < \inf \left\{ \int_\Omega F \d \mu \colon 
 \delta_0 \prec_H \mu , \; \mu \in 
\mc M^+_{\ac} (\Omega)\right\},	
\end{equation}
where $\delta_0$ is the Dirac measure at $0$, then for any  function $\sigma \in C(\Omega)$  
satisfying \eqref{c:rD},
there is a function $h\in H$ such that\footnote{In \cite{Kh011}, we proposed  by definition that the cone $SH(\Omega)$ is not containing the function identically equal to $ -\infty$.} $h\not\equiv -\infty$ on $\Omega$, and
\begin{equation}\label{it:h}
h(z)\leq F^{(\sigma)}(z):=\int_{D(\sigma(z))} F(z+w) \d m^{(\sigma(z))}(w), \quad \forall z\in \Omega.	
\end{equation}
\end{theorem}

The space  $H=\Har (\Omega)$ and the cone $H=SH(\Omega )$  satisfy the
conditions of Theorem \ref{th:5} (see \cite{Kh011}--\cite{Kh012} for example).

Results more general than Theorem \ref{th:5} were proved
for abstract cones $H$ of vector lattices \cite[Theorems 5.1, 6.1]{Kh011}, and, for $H=SH(\Omega )$, in \cite[Corollary 1.7]{CR}.
They require a certain preparation,  and we do not use their here.

We also need another treatments of theorem \ref{th:5} for the cases
 $H=\Har (\Omega )$ and $H=SH(\Omega )$, and for special functions $F$
 representing by the difference \eqref{rep:FuM} 
where $M\in L_{\loc}^1(\Omega)$ or $M\in SH(\Omega)$. 

The following Theorem \ref{cor:termmu} is an earlier version of Theorem \ref{th:5} which joins Proposition \ref{pr:dirmaj} with Theorem \ref{th:5} for special functions \eqref{rep:FuM} with $M\in C(\Omega)$.
\begin{theorem}[{\rm \cite[Main Theorem]{KhIs}}]\label{cor:termmu}
 Let $u\in SH(\Omega  )$ and $M\in C(\Omega )$. The difference \eqref{rep:FuM} 
 admits a harmonic  (subharmonic resp.) minorant on $\Omega $ if and only if
there is a constant $C$ such that 
 $\int_{\Omega} u \d \mu \leq \int_{\Omega} M \d \mu +C$
for all $\mu \in \mc A \mc S (\Omega)$
(for all $\mu \in \mc J (\Omega)$ resp.).
\end{theorem}
We also need an other version of Theorem \ref{cor:termmu}
in the case that the function $M$ in \eqref{rep:FuM}  
belongs to $L^1_{\loc} (\Omega)$ or $SH(\Omega)$.

\begin{theorem}\label{cor:termV}
Let $u \not\equiv -\infty$  be a subharmonic function with the Riesz measure $\nu_u$,   and $u(0)\neq -\infty$. 
Let $M \in L_{\loc}^1(\Omega)$ be a function that is bounded 
on a subdomain $\Omega_1 \Subset \Omega$ containing the  origin. 

If there are  a subdomain $\Omega_0 \Subset \Omega_1$ and 
constant $C$ such that the inequality
\begin{equation}\label{in:arsV}
	\int_{\Omega} u \d \mu \leq \int_{\Omega } M \d \mu +C \quad \text{or}
	\quad
	\int_{\Omega} V_{\mu} \d \nu_u \leq \int_{\Omega } M \d \mu +C 
\end{equation}
holds
for  every $\mu \in \mc A \mc S ^{\Omega_0}(\Omega)\bigcap \mc M_{\ac}(\Omega)$
 (for every   
 $ \mu \in \mc J^{\Omega_0}(\Omega)\bigcap \mc M_{\ac}(\Omega)$ resp.) where $V_{\mu}$ is the potential\/ \eqref{df:poten} of $\mu$, then, for  any function $\sigma \in C(\Omega)$  
 satisfying \eqref{c:rD},
there exists a harmonic (subharmonic resp.) function $h$ 
(\;$h\not\equiv -\infty$ resp.) on\/ $\Omega$
such that (see \eqref{df:si})
\begin{equation}\label{in:hM}
u(z)+h(z)\leq M^{(\sigma)}(z), \quad \forall z\in \Omega .	
\end{equation}
In addition, if $M\in SH(\Omega)$, then the inequality \eqref{in:arsV}
can be replace by 
\begin{equation}\label{in:arsVM}
	\int_{\Omega} V \d \nu_u \leq \int_{\Omega } V \d \nu_M +C 
\end{equation}
for every function $V \in \mc P_{AS}^{\Omega_0}(\Omega)\bigcap C(\Omega)$
 (function $V\in \mc P_{J}^{\Omega_0}(\Omega)\bigcap C(\Omega)$ resp.).
\end{theorem}

\begin{proof} Here and later we need the following elementary
\begin{lemma}\label{lem:o} Let\/ $\Omega_0 \Subset \Omega_1$ 
are subdomains of\/ $\Omega$, and\/ $0\in \Omega_0$. Then there is a domain\/ $\Omega' \in \mc U_0^d(\Omega_1; \Omega_0)$
such that\/ $\Omega_0 \Subset \Omega' \Subset \Omega_1$. In particular,
the domain\/ $\Omega'$ is regular for the Dirichlet problem. 
\end{lemma}
By Lemma \ref{lem:o} we can assume that 
$\Omega_1$ is regular for the Dirichlet problem, $m(\partial \Omega_1)=0$,
and $M$ is bounded in an open neighbourhood of the closure $\overline \Omega_1$. Then we can choose a sufficiently small constant $\e >0$  
and a constant $c \geq 0 $  such that  
\begin{subequations}\label{r:eD}
\begin{align}
|M(z+w)|\leq	c, \quad \forall z\in \Omega_1, \quad 
\forall w\in D(\e), \tag{\ref{r:eD}M}\label{eDM}\\ 
\bigcup_{z\in \Omega_1} D(z, \e) \Subset \Omega ,
\quad 
 \Omega_0 \bigcap \left(\bigcup_{z\in \partial \Omega_1} D(z, \e)\right) =\varnothing .
\tag{\ref{r:eD}$\Omega$}\label{eDMO}
\end{align}
\end{subequations}
We put $F=M-u$, and $H=SH(\Omega)$ or $H=\Har (\Omega)$. 
According to generalized Poisson--Jensen 
formula \eqref{for:PJ} and Proposition \ref{pr:isao}, it follows from \eqref{in:arsV} that the inequality 
\begin{equation}\label{ie:FF}
\int F \d \mu >-C  	
\end{equation}
holds for  every measure $\mu \in \mc M_{\ac}(\Omega)$
such that $\delta_0 \prec_H \mu$ and $\Omega_0 \bigcap \supp \mu 
=\varnothing$. 

Let $\mu$ be an arbitrary measure that belongs to $ \mc M^+_{\ac}(\Omega)$,
and $\delta_0 \prec_H \mu$. Denote by $\mu_1$ the restriction of $\mu$ to $\Omega_1$.
We have the representation
\begin{equation}\label{rep:m}
\mu = \mu_1+\mu_{\infty}, \quad \mu_1, \mu_{\infty} \in  	\mc M^+_{\ac}(\Omega).
\end{equation}
The measure $\mu':=\mu_1^{\Omega_1}*m^{(\e)} \in \mc 
M^+_{\ac}(\Omega)$ is the balayage of $\mu_1$
with respect to $SH(\Omega)$ (see \eqref{eq:AMSJ} and \cite[Lemma 7.1]{Kh012}). 
Hence $\delta_0 \prec_H (\mu'+\mu_{\infty}) \in  \mc M^+_{\ac}(\Omega)$,
\begin{equation}\label{b:u}
\int u \d \mu_1 \leq \int u \d \mu',	
\end{equation}
and, by constructions \eqref{df:baly} and \eqref{eDMO}, $\Omega_0 \bigcap \supp (\mu'+\mu_{\infty}) =\varnothing$. The inequality \eqref{ie:FF}
holds for all such measures. Thus,
\begin{equation}\label{in:FC}
\int F \d (\mu'+\mu_{\infty})\geq -C.	
\end{equation}
Hence we get
\begin{multline*}
\int F \d \mu \overset{\eqref{rep:m}}{=}\int F \d (\mu'+\mu_{\infty})
+\int F \d \mu_1	-\int (M-u) \d \mu' \\
\overset{\eqref{in:FC}}{\geq} -C -\left(\int u \d \mu_1 - \int u \d \mu'\right)
+\int_{\overline \Omega_1} M \d \mu_1- 
\int_{\overline \Omega_1} M \d \mu'\\
\overset{\eqref{b:u}}{\geq}
-C -\int_{\overline \Omega_1} |M| \d (\mu_1 +\mu')
\overset{\eqref{eDM}, \eqref{c:prob}}{\geq}
-C-2c.
\end{multline*}
The last inequality implies the condition  \eqref{i:absb}  
of Theorem \ref{th:5}. 
It follows that there is 
a function $h \in H$ such that \eqref{it:h} holds. Under our notation, \eqref{it:h} consides  with \eqref{in:hM}. 

For the case \eqref{in:arsVM}, by Proposition \ref{pr:isao}, we see that the map $\mc P$ acts from $\mc A \mc S^{\Omega_0}(\Omega)\bigcap \mc M_{\ac}(\Omega)$ into $\mc P_{AS}^{\Omega_0}(\Omega)\bigcap C(\Omega)$
and from $\mc A \mc S^{\Omega_0}(\Omega)\bigcap \mc M_{\ac}(\Omega)$
into $\mc P_{AS}^{\Omega_0}(\Omega)\bigcap C(\Omega)$.
Hence,  \eqref{in:arsV} follows from \eqref{in:arsVM}
according to generalized Poisson--Jensen 
formula \eqref{for:PJ}.   
\end{proof}

\begin{remark} In Theorem \ref{cor:termV} we may assume  that the inequalities \eqref{in:arsV} are fulfilled only for
all measures $\mu \in \mc A \mc S ^{\overline{\Omega}_0}(\Omega)\bigcap \mc M_{\ac}(\Omega)$
 (only for all measures   
 $ \mu \in \mc J^{\overline{\Omega}_0}(\Omega)\bigcap \mc M_{\ac}(\Omega)$ resp.) where
\begin{subequations}\label{df:clO}
\begin{align}
\mc A \mc S ^{\overline{\Omega}_0}(\Omega) &:=\{ \mu \in 
 \mc A \mc S (\Omega) \colon \overline{\Omega}_0 \cap \supp \mu =\varnothing\},  \notag \\  
 \mc J ^{\overline{\Omega}_0}(\Omega) & :=\{ \mu \in 
 \mc J (\Omega) \colon \overline{\Omega}_0 \cap \supp \mu =\varnothing\} .\tag{\ref{df:clO}j}\label{df:clOj}
\end{align}
\end{subequations} 
 It follows from an arbitrary rule for selection of $\Omega_0 \Subset \Omega_1$.  
\end{remark}

\section{From Green functions\\ to Jensen measures}\label{sec:3} 
\setcounter{equation}{0}

The aim of this section is to prove that the conditions of 
statements (Z) and (S) of the Main Theorem give an united condition in terms of Jensen measures.  For the statement (S) this condition is the following evident corollary 
of Proposition \ref{pr:dirmajs}:
\begin{enumerate}
	\item[{\rm (J)}]  
{\it There exists a constant $C$ such that
\begin{equation}\label{c:ihahn}
\int_{\Omega} V_{\mu} \d \nu_u \leq 
\int_\Omega M \d \mu +C 	
\end{equation}
for all measures $\mu \in  \mc J^{\Omega_0}(\Omega)$}.
\end{enumerate}

In subsections \ref{ss:appJ} and \ref{ss:conJ} we show that the condition (J) also is fulfilled under the conditions of (Z).
 
\subsection{A special approximation of Jensen measures}\label{ss:appJ}

First we recall joint results of B.~Cole and T.~Ransford 
\cite{CR01} about the approximation of Jensen measures by harmonic measures.

We write $\conv A$ (\;$\overline{\conv} A$ resp.) for the convex 
(closed convex resp.) hull of a set $A$.

Let $S$ be an open subset or a compact subset of $\Omega$. The space $C(S)$ is a Fr\'echet space with the topology of uniform convergence on  
compact subsets of $S$. The dual space $C(S)^*$ may be identified
with the space of finite signed Borel measures on $\Omega$ of compact support. We use only weak$^*$-topology in the space $C(S)^*$.

\begin{theorem}[{\rm \cite[Theorem 6.6]{CR}}]\label{th:CR}
$\mc  J(\Omega)=\overline{\conv} (\mc H (\Omega))$.
\end{theorem}

\begin{proposition}[{\rm \cite[Proposition 2.1]{CR}, \cite[Proposition 1.1]{Kh03}}]\label{pr:CR1} Let $\Omega_1 $ be a subdomain of $\Omega$, $0\in \Omega_1$. 
Then $\mc J(\Omega_1) \subset \mc J(\Omega)$. If $\mu \in \mc J(\Omega)$ and $\supp \mu \subset \Omega_1$, and  if each bounded component of $\bC \setminus \Omega_1$ meets $\bC \setminus \Omega $, then $\mu \in \mc J(\Omega_1)$.
\end{proposition}

\begin{proposition}\label{pr:appru}
Let\/ $D$ be  a subdomain of\/ $\Omega$, and\/ $0\in D \Subset \Omega$. Then there exists a increasing sequence of domains $D_n \in \mc U_0^d (D)$ such that the sequence of harmonic
measures $\omega_{D_n}(0,\cdot)$ converges to $\omega_D(0,\cdot )$
in $C(\Omega)^*$.
\end{proposition}
\begin{proof} See \cite[Theorems 4.15, 5.14]{Landkof}.
\end{proof}
\begin{proposition}\label{pr:ea}  
Let $\Omega_0$ and $\Omega_1$ are  subdomains of $\Omega$ such that 
\begin{equation}\label{cn:O1}
0\in \Omega_0 \Subset \Omega_1\Subset \Omega	.
\end{equation}
Let $K$ be a compact subset of $\Omega$ satisfying 
$K \bigcap  {\overline{\Omega}_0} =\varnothing$.	
Then  there are domains $ \Omega'\in \mc U_0^d ( \Omega_1; \Omega_0 )$ and $\Omega''\in \mc U_0^d ( \Omega; \Omega_1 )$ such that 
\begin{equation}\label{c:oms}
\Omega_0 \Subset  \Omega' \Subset \Omega_1 \Subset \Omega'' \Subset \Omega , \quad  K \subset \Omega_1\setminus \overline{\Omega'},
\end{equation}
and each bounded component of $\bC \setminus \Omega''$ meets $\bC \setminus \Omega $.
\end{proposition} 
The proof of Proposition \ref{pr:ea} follows from Lemma \ref{lem:o}.

Given domains $\Omega' \Subset \Omega_2 $ containing $0$, we write
\begin{align*}
\mc H^{\ud}( \Omega_2  )&:=
\{ \omega_D(0, \cdot)
\colon D\in \mc U_0^d ( \Omega_2  ) \}, \\
\mc H^{\ud}( \Omega_2 ; \Omega' )&:=\{ \omega_D(0, \cdot)
\colon D\in \mc U_0^d ( \Omega_2 ; \Omega' ) \}.
\end{align*}

\begin{proposition}\label{cn:conm}  
Let $\mu$ be a Jensen measure for $0$ on $\Omega$, and $K:=\supp \mu$.
Under the notations and conditions of Proposition {\rm \ref{pr:ea}},
\begin{equation}\label{it:mu}
\mu \in \overline{\conv} \mc H^{\ud}( \Omega'' ; \Omega' )
\end{equation}
where  the closure is taken with respect to the weak$^*$-topology on $C(\Omega'')^*$. 
\end{proposition} 
\begin{proof} 
By Proposition \ref{pr:ea} (for $\Omega''$ instead of $\Omega$) there exists one more domain $\Omega_2 \in \mc U_0^{d}(\Omega''; \Omega_1 )$ such that $\Omega_1\Subset \Omega_2 \Subset \Omega''$,  
each bounded component of $\bC \setminus \Omega_2$ meets $\bC \setminus \Omega''$, and  $\supp \mu \subset \Omega_2 \setminus \overline{\Omega'}$. 
By Proposition  \ref{pr:CR1} (for $\Omega_2$ instead of $\Omega_1$) 
and Theorem \ref{th:CR} (for $\Omega_2$ instead of $\Omega$),
the measure $\mu$ belongs to $\overline{\conv} \mc H (\Omega_2)$
where the closure is taken with respect to weak$^*$-topology
of $C(\Omega_2)^*$. 
According to  Proposition \ref{pr:appru} we have
\begin{equation*}
\overline{\conv} \mc H (\Omega_2)=
\overline{\conv} \overline{\mc H^{\ud} (\Omega_2)}
=	\overline{\conv} \mc H^{\ud} (\Omega_2)
\end{equation*}
in $C(\Omega_2)^*$. Therefore a net $\{ \sigma_{\gamma}\}
\subset {\conv} \mc H^{\ud} (\Omega_2)$,
$\gamma \in \Gamma$, converges to $\mu$ in $C(\Omega_2)^*$.
For each subscript $\gamma $ there are a finitely many numbers 
$c_k(\gamma)>0$ and domains $D_k(\gamma) \in \mc U_0^d(\Omega_2)$ such that   
\begin{equation}\label{df:init}
\sigma_{\gamma}=
\sum_{k} c_k(\gamma)\omega_{D_k(\gamma)}	, \quad
\omega_{D_k(\gamma)}:=\omega_{D_k(\gamma)}(0, \cdot), \quad  \sum_k c_k(\gamma)=1.
\end{equation}
There is a domain $\Omega_-$ such that $\Omega' \Subset \Omega_-\Subset \Omega_1$ and $K\subset \Omega_2\setminus \overline \Omega_-$.
For every $k, \gamma$ we can choose a domain $\Omega_k(\gamma)\in \mc U_0^{d}(\Omega_-; \Omega' )$ such that the domain $D_k'(\gamma):=D_k(\gamma)\cup \Omega_k(\gamma)\in \mc U_0^d(\Omega_2; \Omega')$ is a regular domain, i.\,e., the complement $\bC \setminus \bigl(D_k(\gamma)\cup \Omega_k(\gamma)\bigr)$ has not of isolated points. Consider the restrictions
\begin{equation}\label{dfg:s}
\omega_{D_k(\gamma)}':=\omega_{D_k(\gamma)}\bigm|_{\Omega_k(\gamma)}, \quad	
\sigma_{\gamma}':=\sum_{k}
c_k(\gamma)\omega_{D_k(\gamma)}'. 
\end{equation}
Since $\supp \mu \Subset \Omega_2\setminus \overline \Omega_-$,
the net $\{ \sigma_{\gamma}'\}$ converges to the null measure in
$C(\overline{\Omega'})^*$:
\begin{equation}\label{con:b'}
\sum_{k} c_k(\gamma)\omega_{D_k(\gamma)}'(\overline{\Omega'})
\underset{\Gamma}{\longrightarrow}0 \, .
\end{equation}
 
Now we reconstruct the net $\{ \sigma_{\gamma}\}$
with the help of balayage.

Denote by 
$(\omega_{D_k(\gamma)}')^{D_k'(\gamma)}$  
the classical balayage of $\omega_{D_k(\gamma)}'$
relative to $D_k'(\gamma)$, and denote by $\omega_{D_k(\gamma)}''$ the classical balayage of $\omega_{D_k(\gamma)}$
relative to $D_k'(\gamma)$. 
By Proposition \ref{pr:harmO} with $D=D_k(\gamma)$ each measure $\omega_{D_k(\gamma)}''$ is exactly the harmonic measure $\omega_{D_k'(\gamma)}(0, \cdot)$. 
So, by definition \eqref{df:baly} and equlities \eqref{eq:harme}, we obtain 
\begin{subequations}\label{dfe:bl}
\begin{align}
&(\omega_{D_k(\gamma)}')^{D_k'(\gamma)}(\Omega_2)=
\omega_{D_k(\gamma)}'(\Omega'), 
\tag{\ref{dfe:bl}}\label{dfe:bl'}\\ 
&\supp (\omega_{D_k(\gamma)}')^{D_k'(\gamma)}\subset {\overline{\Omega}_2}.	
\tag{\ref{dfe:bl}$'$}\label{dfe:bl'a}\\
&\omega_{D_k'(\gamma)}(0, \cdot)=\omega_{D_k(\gamma)}''=(\omega_{D_k(\gamma)}')^{D_k'(\gamma)}
+(\omega_{D_k(\gamma)}-\omega_{D_k(\gamma)}').
\tag{\ref{dfe:bl}$''$}\label{dfe:bl''}
\end{align}
\end{subequations}
It follows from  \eqref{con:b'} and \eqref{dfe:bl'} that
\begin{equation*}
\sum_{k} c_k(\gamma)(\omega_{D_k(\gamma)}')^{D_k'(\gamma)}(\Omega_2) \underset{\Gamma}{\longrightarrow}0 \, . 
\end{equation*}
Hence, in view of \eqref{dfe:bl'a}, \eqref{dfg:s}, and \eqref{con:b'}, the nets 
\begin{equation*}
\sigma_{\gamma}''=\sum_{k} c_k(\gamma)(\omega_{D_k(\gamma)}')^{D_k'(\gamma)},
\quad \sigma_{\gamma}'=\sum_{k} c_k(\gamma)\omega_{D_k(\gamma)}'
\end{equation*}
converge to the null  measure in $C({\overline{\Omega}_2})^*$. It follows from  the representations
\eqref{dfe:bl''} and \eqref{df:init} that the net
\begin{equation*}
\sigma_{\gamma}^*:=\sum_{k} c_k(\gamma)\omega_{D_k'(\gamma)}(0, \cdot)	
=\sigma_{\gamma}''+\sigma_{\gamma}-\sigma_{\gamma}', \quad 
\Omega'\subset D_k'(\gamma)\Subset \Omega_2,
\end{equation*}
converges to $\mu$ in  $C({\overline{\Omega}_2})^*$, and all the more
in $C({{\Omega}''})^*$. So, \eqref{it:mu} is proved.
\end{proof}
 
\subsection{The property (J) for the statement (Z)}\label{ss:conJ}

Suppose that a Jensen measure $\mu$ belongs to $ \mc J^{\overline{\Omega}_0}(\Omega)$ (see 
definition \eqref{df:clOj}) for  the same subdomain $\Omega_0$ as in (Z). Put $K:=\supp \mu$. The function $M$ of Main Theorem is bounded 
on the subdomain $\Omega_1$ of Main Theorem satisfying \eqref{cn:O1}. By Propositions  \ref{pr:ea} and \ref{cn:conm}, there are domains $\Omega', \Omega''\in \mc U_0^d ( \Omega_1; \Omega_0 )$  such that \eqref{c:oms} and \eqref{it:mu} are fulfilled. 

Put 
\begin{equation}\label{ch:dd}
d':=\dist (\Omega_0, \partial \Omega' ), \quad 
d'':=\dist (\Omega'', \partial \Omega ).	
\end{equation}
We can choose a number $\e$ so that
\begin{equation}\label{ch:ee}
0< \e < \min\{ d', d''\}.	
\end{equation}
Let $D \in \mc U_0^d( \Omega'';\Omega')$.  According to
\eqref{ch:dd} and \eqref{ch:ee}, for any $w \in D(\e )$, 
the domain
$D_{w}=\{ z-w \colon  z\in D\}$ 	
belongs to the class $\mc U_0^d( \Omega;\Omega_0)$. Therefore,
by the main condition of (Z), there is a constant  $C$
 such that the inequality
\begin{equation}\label{in:masn}
\int_{\Omega} g_{D_{w}}(\zeta , 0) \d \nu_u(\zeta)\leq
\int_{\Omega} g_{D_{w}} (\zeta , 0 )\d \nu_M (\zeta )+C 
\end{equation}
holds for each domain 
$D \in \mc U_0^d ( \Omega'' ; \Omega' )$.
Here the constant $C$ is independent of $D$, $w$, and $\mu$.
 
By the generalized Poisson--Jensen formula \eqref{for:PJ},
the inequality \eqref{in:masn} can be rewritten as 
\begin{equation}\label{in:masnn}
\int_{\Omega} u(z)\d \omega_{D_{w}}(0 , z) \leq
\int_{\Omega} M(z)\d \omega_{D_{w}} (0, z )+C' 
\end{equation}
where the constant $C'= C-u(0)-M(0)$ is independent of 
$D\in \mc U_0^d( \Omega'';\Omega')$, $w\in D(\e)$, $\mu$, and $\e$ under condition \eqref{ch:ee}.
If we replace $z-w$ by a new variable in \eqref{in:masnn}, then
\begin{equation}\label{in:maonn}
\int_{\Omega} u(z+w)\d \omega_{D}(0 , z) \leq
\int_{\Omega} M(z+w)\d \omega_{D} (0, z )+C' 
\end{equation} 
for all $D \in  \mc U_0^d( \Omega'';\Omega')$ and $w\in D(\e)$
where the constant $C'$ is independent of $\mu$ and $\e$ of  \eqref{ch:ee}.
Integrating the last inequality with respect to the probability measure $m^{(\e)}$ (see \eqref{df:mr}), by Fubini's theorem, we obtain
\begin{multline*}
 \int_{\Omega} \bigl(u*m^{(\e)}\bigr)(z)\d \omega_{D}(0 , z) 
=\int_{D(\e)} \int_{\Omega} u(z+w)\d \omega_{D}(0 , z) \d m^{(\e)}(w)
\\ \overset{\eqref{in:maonn}}{\leq}
\int_{D(\e)}\int_{\Omega} M(z+w)\d \omega_{D} (0, z ) \d m^{(\e)}(w)+
\int_{D(\e)}C' \d m^{(\e)}(w)\\=
\int_{\Omega} \bigl(M*m^{(\e)}\bigr)(z)\d \omega_{D} (0, z ) +
C' . 
\end{multline*}
Hence, under notation 
\eqref{df:si}, 
\begin{equation*}
\int_{\Omega} u^{(\e)}(z)\d \omega_{D}(0 , z)
\leq
\int_{\Omega} M^{(\e)}(z)\d \omega_{D} (0, z )+C' 
\end{equation*}
for all domains $D\in \mc U_0^d( \Omega'';\Omega')$ where the functions  $u^{(\e)}$ and $M^{(\e)}$ are well defined and continuous on $\Omega''$. According to \eqref{it:mu} of Proposition \ref{cn:conm}, 
it implies that  
\begin{equation*}
\int_{\Omega} u^{(\e)}(z)\d \mu( z)
\leq
\int_{\Omega} M^{(\e)}(z)\d \mu ( z )+C'. 
\end{equation*}
In view of \eqref{df:convo}, we have
\begin{equation*}
\int_{\Omega} u\d \mu
\leq
\int_{\Omega} M^{(\e)}\d \mu +C' . 
\end{equation*}
If $\e$ tend to $0$, then the decreasing net $\{ M^{(\e)}\}$ of continuous functions tend pointwise to $M$ on $\Omega''$ since $M$ is subharmonic. Thus, 
\begin{equation}\label{in:itogo}
\int_{\Omega} u\d \mu
\leq
\int_{\Omega} M \d \mu +C'. 
\end{equation}
  
But the measure $\mu \in \mc J^{\overline{\Omega}_0}(\Omega)$ was any given, and the constant $C'$ is independent of $\mu$. Therefore, the inequality \eqref{in:itogo} is fulfilled
for every measure $\mu \in \mc J^{\overline{\Omega}_0}(\Omega)$. By the generalized   Poisson--Jensen formula \eqref{for:PJ}, the inequality 
\eqref{in:itogo} implies the inequality 
\begin{equation}\label{in:itogi}
\int_{\Omega} V_{\mu} \d \nu_u
\leq
\int_{\Omega} M \d \mu +\bigl(C'-u(0)\bigr)  
\end{equation}
for every measure $\mu \in \mc J^{\overline{\Omega}_0}(\Omega)$ where $V_{\mu}$ is potential \eqref{df:poten}
of $\mu$, and the constant $C'':=C'-u(0)$ is independent of $\mu$. 
But the inequality \eqref{in:itogi} is exactly
\eqref{c:ihahn} with constant $C''$ instead of $C$. Besides, here we can replace the uper index  of $\mc J^{\overline{\Omega}_0}(\Omega)$ by  $\Omega_0$. The last follows from an arbitrary rule for selection of $\Omega_0 \Subset \Omega_1$. The property (J) is proved.$\,\bullet$ 

\begin{remark} For the deduction of (J) from (S) or (Z)
we did not use a some representation of $M$ by a subharmonic kernel. 
In other words, the conditions  \ref{rmii}) and \ref{rmiii}) 
of Main Theorem are unnecessary in order to prove (J).	   
\end{remark}

\section{The proofs of the Main Theorem}\label{sec:4} 
\setcounter{equation}{0}

\subsection{From Jensen  measures and functions\\ to
the existence of harmonic minorant}

The aim of this subsection is to prove that the condition (J)
at the beginning of Section \ref{sec:3} implies the existence of harmonic minorant for the function \eqref{rep:MQu}. After Section \ref{sec:3} it gives both (S) and  (Z). 

Suppose that the condition (J) is fulfilled. 

By the generalized Poisson--Jensen formula \eqref{for:PJ}
and Proposition \ref{pr:isao}, it follows from the condition (J) that, for the constant $C'=C+M(0)$, 
\begin{enumerate}
	\item[($\mathrm J'$)] {\it the inequality 
\begin{equation}\label{c:ihahJ}
\int_{\Omega } V \d \nu_u \leq 
\int_\Omega  V \d \nu_M +C' 	
\end{equation}
holds for every {\it Jensen function\/} $V\in \mc P_{J}^{\Omega_0}(\Omega)$} (see above \eqref{li:V}). 
\end{enumerate}

Now let $\mu \in \mc A \mc S^{\Omega_0} (\Omega)\cap \mc M^+_{\ac}(\Omega )$ be a {\it Arens--Singer measure with the potential $V_{\mu}$} 
(see \eqref{df:poten}). By Proposition \ref{pr:isao} the Arens--Singer function $V_{\mu}$ belongs to the class $\mc P_{AS}^{\Omega_0} (\Omega)$. Besides, by Definition  \ref{df:asjf}, the 
function 
\begin{equation}
V_{\mu}^+(\zeta ):=\max \{ V(\zeta), 0\}, \quad \zeta \in \bC \setminus \{ 0\},	
\end{equation}
is a {\it Jensen function\/} that belongs to the class $\mc P_{J}^{\Omega_0} (\Omega)$. Therefore the condition ($\mathrm J'$) implies that the inequalities
\begin{equation}\label{est:V+}
\int_{\Omega } V_{\mu} \d \nu_u \leq 
\int_{\Omega } V_{\mu}^+ \d \nu_u
\leq \int_\Omega  V_{\mu}^+\d \nu_M +C'	
\end{equation}
are fulfilled {\it for all  Arens--Singer measures  
$\mu \in \mc A \mc S^{\Omega_0} (\Omega)$}.

We shall need to estimate above the last integral 
\begin{equation}\label{df:Im}
I_{\mu}:=\int_\Omega  V_{\mu}^+\d \nu_M\overset{\rm \ref{rmii})}{=}\int_B  V_{\mu}^+\d \nu_M.	
\end{equation}

According to the representation  \eqref{repr:Vmu} of Proposition \ref{prop:Vmu} for the potential $V_{\mu}$ we have
\begin{multline*}
	I_{\mu}=\int_{B}\left(\int_{\Omega}
 \bigl( k(\zeta , z)-k(\zeta , 0)\bigr)  \d \mu (z)\right)^+
\d \nu_M (\zeta )\\
\leq \int_{B}\int_{\Omega}
 \bigl( k(\zeta , z)-k(\zeta , 0)\bigr)^+  \d \mu (z)
\d \nu_M (\zeta )\\
=\int_{B}\int_{\Omega}
 \left(\bigl( k(\zeta , z)-k(\zeta , 0)\bigr) 
+
 \bigl( k(\zeta , 0)-k(\zeta , z)\bigr)^+ \right)\d \mu (z)
 \d \nu_M (\zeta )
\end{multline*}
Hence, by Fubini's theorem, 
\begin{multline}\label{for:ineq}
I_{\mu}\leq 
 \int_{\Omega}
 \left(  \int_{B}\bigl( k(\zeta , z)-k(\zeta , 0 )\bigr) 
\d \nu_M (\zeta )\right.\\
\left.+
\int_{B} \bigl( k(\zeta , 0)-k(\zeta , z)\bigr)^+ 
 \d \nu_M (\zeta )\right)\d \mu (z)\\
\overset{\eqref{df:Qkn}}{=}\int_{\Omega}
 \left(  \int_{B}\bigl( k(\zeta , z)-k(\zeta , 0 )\bigr) 
\d \nu_M (\zeta ) +Q^{\nu_M}_k (z) \right)\d \mu (z).
\end{multline}
By condition \ref{rmii}) and representation \eqref{repr:Riesz} of Proposition \ref{pr:Riesz} the function $M$ can be represented in the form
\begin{equation*}
M(z)=\int_{B} k(\zeta , z) \d \nu_M(\zeta )+H(z)
=U_k^{\nu_M}(z)+H(z) , \quad z\in \Omega ,	
\end{equation*}
where $H \in \Har (\Omega)$. From this it follows that
\begin{equation}\label{lh:s}
\bigl(M(z)-M(0)\bigr)-\bigl(H(z)-H(0)\bigr)=	\int_{B}\bigl( k(\zeta , z)-k(\zeta , 0 )\bigr) 
\d \nu_M (\zeta ).
\end{equation}
Substituting the left-hand side of \eqref{lh:s} in \eqref{for:ineq}, we get 
\begin{equation}\label{sd:ll}
I_{\mu} \overset{\rm \ref{rmiii})}{\leq} \int_{\Omega}
\Bigl(\bigl(M(z)-M(0)\bigr)-\bigl(H(z)-H(0)\bigr)+Q(z)	\Bigr)\d \mu (z).
\end{equation}
By Definition \ref{df:asj} of Arens--Singer (probability \eqref{c:prob}) measures the right-hand side
of \eqref{sd:ll} is equal to
\begin{multline*}
	\int_{\Omega}
\Bigl(\bigl(M(z)+Q(z)\bigr)-H(z)\Bigr)\d \mu (z)
-M(0)+H(0)\\
\overset{\eqref{df:bal1}}{=} 
\int_{\Omega} \bigl(M(z)+Q(z)\bigr) \d \mu (z) -M(0)
\end{multline*}
since $H\in \Har (\Omega)$. Now if we recall \eqref{est:V+}, \eqref{df:Im}, and \eqref{sd:ll}, we get
\begin{equation*}
\int_{\Omega } V_{\mu} \d \nu_u \leq 
\int_{\Omega} \bigl(M(z)+Q(z)\bigr) \d \mu (z) +( C'-M(0))
\end{equation*}
where the constant $C'':=C'-M(0)$ is independent of $\mu \in \mc A \mc S^{\Omega_0} (\Omega)\cap \mc M^+_{\ac}(\Omega )$.
It means that second inequality in \eqref{in:arsV} is fulfilled {\it for every  Arens--Singer measure\/} $\mu \in \mc A \mc S^{\Omega_0} (\Omega)\cap \mc M^+_{\ac}(\Omega )$ with the the function $M+Q \in L_{\loc}^1(\Omega)$ instead of $M$, and with the  constant $C''$ instead of $C$. Thus, by Theorem \ref{cor:termV},   
for  any function $\sigma \in C(\Omega)$  satisfying \eqref{c:rD},
there exists a function $h \in \Har (\Omega )$ such that 
\begin{equation*}
u(z)+h(z)\leq (M+Q)^{(\sigma)}(z)=\bigl(M^{(\sigma)}+Q^{(\sigma)}\bigr)(z), \quad \forall z\in \Omega .	
\end{equation*}
In other words, the function $h$ is a harmonic minorant 
for the function \eqref{rep:MQu}.

This completes the proof of the Main Theorem. $\,\bullet$

\subsection{The proof of Corollary \ref{th:Mc}}

In this subsection we use a some function $\sigma \in C(\Omega)$
satisfying  \eqref{c:rD}.  
\subparagraph{The proof of {\rm ($\mathrm Z_{\sic}$).}}
Let $f_\Lambda$  be a holomorphic function with zero set 
$\Lambda =\{ \lambda_n\}$, which exists by Weierstrass' theorem.
Then the function $u=\log |f_{\Lambda} |$ is  subharmonic with the Riesz measure $n_{\Lambda}$ where $n_{\Lambda}$ is the counting measure of sequence $\Lambda$ (see \eqref{df:nLa}).
At that rate the inequality \eqref{in:maz} is exactly
\eqref{in:ma} with the Riesz measure $\nu_u=n_{\Lambda}$ of $u$. By assertion  {\rm (Z)} of the Main Theorem there exists a function $h\in \Har (\Omega)$ such that the inequality \eqref{e:estM} holds. The domain $\Omega$ is simply connected. Therefore there exists a holomorphic function $g$ on $\Omega$ such that $h=\Re g$ on $\Omega$. Thus there is a holomorphic function $g$ such that  
\begin{equation}\label{e:estMf}
\log |f_{\Lambda}(z)|+\Re g(z)=u(z)+h(z)
\overset{\eqref{e:estM}}{\leq}  M^{(\sigma)}(z)+Q^{(\sigma)}(z), \quad \forall z\in \Omega. 	
\end{equation}
In other words the function $f=f_{\Lambda}\exp g$
with ${\Zero}_f=\Zero_{f_\Lambda}=\Lambda$
satisfies the inequality $\log |f| \leq M^{(\sigma)}+Q^{(\sigma)}$
on $\Omega$. This proves that $\Lambda$ is a zero set for  
$\Hol(\Omega; M^{(\sigma)}+Q^{(\sigma)})$. $\bullet$
\subparagraph{The proof of {\rm ($\mathrm S_{\sic}$).}}
If $\Lambda$ is zero subset for $\Hol (\Omega; M)$, then  
there exists a holomorphic function $f\not\equiv 0$ 
on $\Omega$ such that $f(\Lambda )=0$, and
$f\in \Hol (\Omega; M)$. 
Let $f_\Lambda$  be a holomorphic function with zero set 
$\Lambda$. Evidently, we have a representation $f=f_{\Lambda}g$ where $g\not\equiv 0$ is a holomorphic function on $\Omega$. 
The condition $f\in \Hol (\Omega; M)$ implies an estimate
\begin{equation}\label{est:enits}
\log |f_{\Lambda}(z)|+\log |g(z)|\leq
M(z)+C, \quad \forall z\in \Omega ,	
\end{equation}
where $C$ is a constant.  We can rewrite \eqref{est:enits}
as 
\begin{equation*}
h(z):=\log |g(z)|-C\leq
M(z)-\log |f_{\Lambda}(z)|, \quad \forall z\in \Omega .	
\end{equation*}
In particular, it means that the difference $M-\log |f_{\Lambda}|$
admits a subharmonic minorant on $\Omega$. Therefore, by statement (S) 
of the Main Theorem (with $\log |f_{\Lambda}|$ instead of $u$), the function 
$M^{(\sigma)}+Q^{(\sigma)}-\log |f_{\Lambda}| $
admits a harmonic minorant $h$ on $\Omega$.
Since the domain $\Omega$ is simply connected, 
there exists a holomorphic function $g$ on $\Omega$ such that $h=\Re g$ on $\Omega$. Thus we get \eqref{e:estMf}. 
As before, this shows that $\Lambda$ is a zero set for  
$\Hol(\Omega; M^{(\sigma)}+Q^{(\sigma)})$. $\bullet$

\subparagraph{The proof of {\rm ($\mathrm M_{\sic}$).}}
Let $f=g/q$ be a meromorphic function on $\Omega$. 
If $g,q \in \Hol (\Omega ; M )$, then the subharmonic function
$u:=\max \bigl\{ \log |g|, \log |q|\bigr\}$ satisfies
\begin{equation}\label{ess:fou}
\max \bigl\{ \log |g(z)|, \log |q(z)|\bigr\}=u(z)\leq M(z)+C, 
\quad \forall z\in \Omega,	
\end{equation}
where $C$ is a constant. It is clear that there is
a representation 
\begin{equation}\label{reply:f}
f=\frac{g_1}{q_1}=\frac{g_1l}{q_1l} \, , \quad g=g_1l , \; q=q_1l, \quad l\in \Hol (\Omega),	\; l\not\equiv 0,
\end{equation}
where $g_1, q_1$ are holomorphic functions without common zeros on $\Omega$.
Hence, in view of \eqref{ess:fou}, we have
\begin{equation}\label{ess:foua}
\max \bigl\{ \log |g_1(z)|, \log |q_1(z)|\bigr\}
+\log |l(z)|=u(z)\leq M(z)+C	
\end{equation}
for all $z\in \Omega$. Put 
\begin{equation}\label{ddf:u1}
u_1:=\max \bigl\{ \log |g_1|, \log |q_1|\bigr\}. 
\end{equation}
Under this notation we can rewrite \eqref{ess:foua} as 
$\log |l| -C \leq M-u_1$ on $\Omega$.
It means that the difference $M-u_1$ admits the subharmonic minorant 
$\log |l|-C$. Therefore, by assertion (S) of the Main Theorem,
the function 
$M^{(\sigma)}+Q^{(\sigma)}-u_1 $ admits a harmonic minorant.
In other words there exists a holomorphic function $s$
such that $u_1+\Re s \leq M^{(\sigma)}+Q^{(\sigma)}$.  
Thus according to  \eqref{ddf:u1},
\begin{equation*}
\max \bigl\{ \log |g_1e^s|, \log |q_1e^s|\bigr\}\leq M^{(\sigma)}+Q^{(\sigma)}.	
\end{equation*}
Hence if we put $g_0=g_1e^s$ and $q_0=q_1e^s$, then in view of 
\eqref{reply:f} $f=g_0/q_0$, $g_0, q_0 \in \Hol (\Omega;M^{(\sigma)}+Q^{(\sigma)})$, and $\Zero_{g_0}\cap \Zero_{q_0}=\varnothing$
as desired.$\,\bullet$

\section{The proof of Theorem \ref{th:MB}}\label{sec:5} 
\setcounter{equation}{0}

\subsection{The reduction\\ to an upper estimate
of function \eqref{df:Qkn}}\label{subsec:alpha}
Since $M$ is bounded below at the origin, we have
$$
\int_0^{1/2} \frac{\nu_M^{\rad} (t)}{t} \d t <+\infty , 
\quad \nu_M(\{0\})=0.
$$
It and the condition \eqref{cond:mea}
imply that the subharmonic Bomash's kernel
  \eqref{df:HB} from the item (${\overline{ \mathrm B}}_s$) of Section
\ref{sec:1}  is suitable for the Riesz measure $\nu_M$ 
of $M$ (see \eqref{c:Bok}).    
We use the subharmonic Bomash's kernel
$\overline b_2(\zeta , z)$ in the role of subharmonic kernel $k$ 
supported by $\D \setminus \{ 0\}$. In this case three conditions 
{\rm \ref{rmi})--\ref{rmiii})} of the Main Theorem are fulfilled.  
It follows from Corollary \ref{th:Mc} that all statements
($\mathrm Z_{\D}$), ($\mathrm S_{\D}$), and ($\mathrm M_{\D}$)
of Theorem  \ref{th:MB} are fulfilled for the space
$\Hol (\D; M^{(\sigma)}+Q^{(\sigma)})$ where
\begin{equation}\label{form:Qb2}
Q(z):=\int_{\D} \bigl(
\overline b_2(\zeta , 0)-\overline b_2(\zeta , z)\bigr)^+  \d {\nu}_M (\zeta ), \quad z\in \Omega.
\end{equation}
By construction this function $Q$ belongs to $L^1_{\loc}(\D)$ since the function $M$ is bounded below on every $D \Subset \D$.
\begin{lemma}\label{lem:Qest} There is an absolute constant $a$
such that, for any\/ $0<\e <1$, 
\begin{equation}\label{es:Qb2}
Q(z)\leq A_{M, \e}(z)-M(z)+C_{\e}\cdot b_M^{[\alpha']}(z), \quad  a\leq |z|<1 ,	
\end{equation}
where the functions\/ $A_{M, \e}$ and\/  $b_M^{[\alpha]}$ is defined by\/ 
\eqref{df:AMe} and\/ \eqref{df:bMe}, $\alpha' >0$ is an absolute constant,
and $C_{\e}$ depends only on $\e$.
\end{lemma}
We shall prove the important upper estimate \eqref{es:Qb2} later by a few steps.

Now we choose a funcion $\sigma \in C(\D )$ satisfying \eqref{c:rD}
so that 
\begin{equation}\label{ch:ev}
\sigma (z)\leq \e (1-|z|), \quad
\bigl(A_{M, \e}\bigr)^{(\sigma)}(z)\leq A_{M, \e}(z)+1,
 \quad \forall z\in \D . 	
\end{equation}
It is possible since the function $A_{M, \e}$ is continuous on $\D$
if $M$  is subharmonic.
It follows from \eqref{ch:ev} and \eqref{es:Qb2} 
that 
\begin{multline*}
M^{(\sigma)}(z)+Q^{(\sigma)}(z)	
\leq  A_{M, \e}(z)+ \bigl(A_{M, \e}(z)-M(z)+C_{\e}\cdot b_M^{[\alpha']}(z)\bigr)^{(\sigma)}
\\ \leq A_{M, \e}(z)+\bigl(A_{M, \e}\bigr)^{(\sigma)}(z) -M^{(\sigma)}(z)+C_{\e}\cdot \bigl(b_M^{[\alpha']}\bigr)^{(\sigma)}(z) 
\\ \leq
2A_{M, \e}(z)-M(z) +C_{\e}\cdot b_M^{[\alpha'+1]}(z)+1, \quad a<|z|<1.
\end{multline*}
Besides, the left-hand side is bounded below in $D(a)$. Thus the statements 
($\mathrm Z_{\D}$), ($\mathrm S_{\D}$), ($\mathrm M_{\D}$)
hold for the space \eqref{sp:defq} with absolute constant $\alpha =\alpha'+1$.$\;\bullet$ 

\subsection{The proof of Lemma \ref{lem:Qest}}

It is enough to estimate of integral  \eqref{form:Qb2} 
in the case $z=x\in [0, 1)$. First we investigate sizes of sets
(see Examples, \eqref{n:bbr} and ($\mathrm B_0$), (${\overline{ \mathrm B}}_s$) and \eqref{df:HB})
\begin{subequations}\label{df:D12}
\begin{align}
D_1(x)&:=\bigl\{ \zeta \in \D \colon 
\overline b_1(\zeta , 0)> \overline b_1(\zeta , x)\bigr\}=
\bigl\{ \zeta \in \D \colon \rho (\zeta , x) < |\zeta| \bigr\},
\tag{\ref{df:D12}.1}\label{df:D121}\\
D_2(x)&:=\bigl\{ \zeta \in \D \colon 
\overline b_2(\zeta , 0)> \overline b_2(\zeta , x)\bigr\}
\notag \\
&\overset{\eqref{df:HBb}}{=}\left\{ \zeta \in \D \colon 
|\zeta|^2(2-|\zeta|^2) > \bigl| \overline B_\zeta ( x)\bigr| \bigl|2- \overline B_\zeta ( x)\bigr|\right\}.
\tag{\ref{df:D12}.2}\label{df:D122}
\end{align}
\end{subequations}
The domain $D_2(x)$ determines the function $Q$ from \eqref{form:Qb2} since
\begin{equation}\label{form:Qb2x}
Q(x)=\int_{D_2(x)} \bigl(
\overline b_2(\zeta , 0)-\overline b_2(\zeta , x)\bigr)  \d {\nu}_M (\zeta ), \quad x\in [0, 1).
\end{equation}
Besides, the function $t \rightarrow t^2(2-t^2)$ is strictly increasing on $[0, 1)$, and $\bigl| \overline B_\zeta ( x)\bigr|<1$ for $\zeta \in \D$.
Therefore, if $\zeta \in D_2(x)$,  then, by \eqref{df:D122}, we have $|\zeta|^2 > \bigl| \overline B_\zeta ( x)\bigr|$, i.\,e., $|\zeta| > \rho (\zeta , x)$. Thus,
\begin{equation}\label{inc:D12}
D_2(x)\subset D_1(x), \quad \forall x\in [0,1).	
\end{equation}

\subparagraph{The domain $D_1(x)$.} We put
\begin{equation}\label{df:zetaa}
	\zeta :=te^{i\theta}, \quad 0\leq t <1, \quad \theta \in (-\pi , \pi ] .
\end{equation}
If $\rho (\zeta , x) \leq |\zeta|$, then 
\begin{equation*}
\Bigl| \frac{te^{i\theta}-x}{1-txe^{-i\theta}}\Bigr|^2=
\frac{t^2-2tx\cos \theta +x^2}{1-2tx\cos \theta +t^2x^2}< t^2. 	
\end{equation*}
The last inequality is equivalent to 
\begin{equation}\label{in:fth}
\cos \theta > \frac{x(1+t^2)}{2t} \, .	
\end{equation}
Hence $\cos \theta >x$ and
\begin{equation}\label{con:fortet}
|\theta| < \frac{\pi}{2} \, \sqrt{1-x^2} .	
\end{equation}

Besides, inequality \eqref{in:fth} implies
$x(1+t^2)< 2t$. 
Hence, for $t=|\zeta|$ of \eqref{df:zetaa}, we obtain
\begin{equation}\label{con:fort}
1>t> \frac{1-\sqrt{1-x^2}}{x} \geq x-\sqrt{1-x^2}, \quad 0 \leq x <1. 	
\end{equation}
The inequalities \eqref{con:fortet} and \eqref{con:fort} give 
the inclusion
\begin{equation}\label{m:D1}
D_1(x)\subset \left\{ \zeta =te^{i\theta} \colon x-\sqrt{1-x^2}<t<1, \;
|\theta| < \frac{\pi}{2} \, \sqrt{1-x^2} \right\}	.
\end{equation}
But this inclusion is insufficient for  good estimates of integral \eqref{form:Qb2x} by \eqref{inc:D12}. It is important farther what, according to \eqref{inc:D12} and \eqref{m:D1},  the set $D_2(x)$ lies {\it in the right half-plane\/} for $x>0$. 
 
\subparagraph{The domain $D_2(x)$.} Let $\zeta \in D_2(x)$ the same as in \eqref{df:zetaa}. 
\begin{proposition}\label{prop:D2x} Under the condition\/ $9/10 \leq x<1$ the Carleson box 
of relative size\/ $6$ with center at\/ $x$ (\,see \eqref{df:box}\,)
\begin{equation}\label{est:D2x}
 {\Box}_6(x)=\bigl\{ \zeta=te^{i\theta} \colon x-6(1-x)<t<1, \; |\theta|< 6(1-x)\bigr\}	
\end{equation}
includes the domain\/ $D_2(x)$.
\end{proposition}
\begin{proof} In view of \eqref{df:HBb}, we have 
$$
\log \bigl(|\zeta|(2-|\zeta|^2) \bigr) >\log \frac{|\zeta -x|  \bigl|2-|\zeta|^2-\overline{\zeta}x\bigr|}{|1-\overline{\zeta}x |^2}
$$
for $\zeta=te^{i\theta} \in D_2(x)$. Hence, under the notation \eqref{df:zetaa},
\begin{multline}\label{in:lonsi}
t^2(2-t^2)^2>
\frac{(t^2-2tx \cos \theta +x^2)\bigl((2-t^2)^2-2(2-t^2)tx \cos \theta +t^2x^2\bigr)}{(1-2tx\cos \theta +t^2x^2)^2}\\
\! =\frac{\bigl((t-x)^2+4tx \sin^2(\theta/2)\bigr)
\bigl((2-t^2-tx)^2+4(2-t^2)tx \sin^2(\theta/2)\bigr)}
{\bigl((1-tx)^2+4tx \sin^2(\theta/2)\bigr)^2} \, .
\end{multline}
In particular, $t>0$. For convenience, we put
\begin{equation}\label{df:sssi}
	s:=\sin^2 \frac{\theta}{2} \, .
\end{equation}
Besides, $(2-t^2)\geq t^2(2-t^2)^2$ for $t\in [0, 1)$. Therefore,
the inequality \eqref{in:lonsi} implies
\begin{equation*}
2-t^2>\frac{\bigl((t-x)^2+4tx s\bigr)
\bigl((2-t^2-tx)^2+4(2-t^2)tx s\bigr)}
{\bigl((1-tx)^2+4tx s\bigr)^2} \, ,	
\end{equation*}
whence a direct calculation gives
\begin{multline}
s\cdot q(t,x):=s \cdot 4tx\bigl((2-t^2-tx)^2+(t-x)^2(2-t^2)-2(1-tx)^2(2-t^2)\bigr)
\\ <(1-tx)^4(2-t^2)-(t-x)^2(2-t^2-tx)^2=:p(t,x),	
\end{multline}
where the polynomials $q$ and $p$ are defined by first and second equalities respectively. We have  the following factorial expansion
for co-factor $q$ after $s$:
\begin{equation*}
q(t,x)=8x^3t(1-t^2)^2>0, \quad t\in (0, 1), \quad x\in (0, 1).	
\end{equation*}
It permits to  get an upper estimate for $s$:
\begin{multline}\label{est:sgx}
s< \frac{p(t,x)}{q(t,x)} =\frac{(1-tx)^4(2-t^2)-(t-x)^2(2-t^2-tx)^2}{8tx^3(1-t^2)^2} \\
=\frac{(1-t)^2(1+t^2)(-t^2-x^4t^2+4x^3t+2-4x^2)}{8x^3t(1-t^2)^2} \\
=\frac1{8x^3t} \, (-t^2-x^4t^2+4x^3t+2-4x^2)\\
=\frac1{8x^3t}\, \bigl( -(1+x^4)t^2+(4x^3)t+2(1-2x^2)\bigr)=:
\frac1{8x^3t}\, g_x(t) \, ,
\end{multline}
where the  quadratic trinomial $g_x$ is defined by last equality. 
By definition \eqref{df:sssi}, $s \geq 0$. Therefore 
the quadratic trinomial $g_x(t)$ is strictly positive
for $\zeta =te^{i\theta}\in D_2(x)$, whence
\begin{multline}\label{es:ttxa}
t>\frac{2x^3-\sqrt{2}(1-x^2)}{1+x^4}=
x-\left(x-\frac{2x^3-\sqrt{2}(1-x^2)}{1+x^4}\right)
\\
=x-(1-x^2)\, \frac{\sqrt{2}+x(1-x^2)}{1+x^4} \, \geq
x-2(1-x^2)\geq x-4(1-x). 
\end{multline}
In particular, 
\begin{equation}\label{est:tfx}
\text{if $x\geq 9/10$, \; then $t > 1/2$}.
\end{equation}

  Now we must  find an upper estimate for $s$.  
Let us consider the inequality \eqref{est:sgx} again. 
The  quadratic trinomial $g_x$ attains its maximum at the point 
$t_x=2x^3/(1+x^4)$ so that 
\begin{multline*}
g_x(t_x)=-(1+x^4)\left(\frac{2x^3}{1+x^4}\right)^2+
4x^3 \, \frac{2x^3}{1+x^4}+2(1-2x^2)\\
=\frac{2-4x^2+2x^4}{1+x^4}=\frac{2(1-x^2)^2}{1+x^4}\leq 
2(1-x)^2\max_{0\leq x <1}\frac{(1+x)^2}{1+x^4}\leq
8(1-x)^2 
\end{multline*}
for each $x\in [0, 1)$. Hence, in view of \eqref{df:sssi} and \eqref{est:sgx}, we get
\begin{equation*}
	\sin^2 \frac{\theta}{2} \, = s < \frac1{8x^3t}\cdot 8(1-x)^2
=\frac1{x^3t}\cdot (1-x)^2
\end{equation*}
whence, for $x\geq 9/10$, by \eqref{est:tfx} we obtain 
\begin{equation}
|\theta| \leq \pi \sqrt{\frac{10^3}{9^3 \cdot (1/2)}} \, (1-x) 
<6(1-x).
\end{equation}
Thus, taking into account \eqref{es:ttxa}, we have the inclusion 
\begin{equation}
D_2(x) \subset \bigl\{ \zeta=te^{i\theta} \colon x-4(1-x)<t<1, \; |\theta|< 6(1-x)\bigr\}	
\end{equation}
for all $9/10 \leq x<1$. This completes the proof of Proposition \ref{prop:D2x}. 
\end{proof}

\subparagraph{The estimates of kernel of integral \eqref{form:Qb2x}.}
Let $0<\e <1$.
We put 
\begin{equation}\label{df:delex}
	\Delta_{\e}(x):=D(x, \e (1-x)).
\end{equation}
Denote by 
\begin{equation}\label{df:gre}
g_{\Delta_{\e}(x)}(\zeta , x):=\log \Bigl| \frac{\e (1-x)}{\zeta -x}\Bigr|
\end{equation}
 the Green's function for the disk $\Delta_{\e}(x)$ with the pole $x$.
\begin{proposition}\label{pr:dexes} If\/      $9/10 \leq x<1$, then  the inequality
\begin{equation}\label{est:b+}
\bigl(\overline b_2(\zeta , 0)-\overline b_2(\zeta , x)\bigr)^+
\leq 	g_{\Delta_{\e}(x)}(\zeta , x)+\log \frac{30}{\e}
\end{equation}
holds for all $\zeta \in \Delta_{\e}(x)$. Besides,
\begin{equation}\label{est:fr2}
\frac{(1-|\zeta|)^2}{(1-x)^2} \geq 1-\e, \quad \forall \zeta \in \Delta_{\e}(x).	
\end{equation}
\end{proposition}  
\begin{proof} By \eqref{df:HBb} and \eqref{df:gre} we obtain
the following representation
\begin{equation}\label{df:nrep}
\overline b_2(\zeta , 0)-\overline b_2(\zeta , x)=
g_{\Delta_{\e}(x)}(\zeta , x)+\log\frac{|\zeta|\bigl(2-|\zeta|^2\bigr)|1-{\overline \zeta}x|^2}{\e (1-x)\bigl|2-|\zeta|^2-\overline{\zeta}x\bigr|}\, .	
\end{equation}
Given $\zeta \in \D$ and $0\leq x <1$, we have $|\zeta|\bigl(2-|\zeta|^2\bigr)\leq 2$,
\begin{equation}\label{est:2zet}
\bigl|2-|\zeta|^2-\overline{\zeta}x\bigr|=
\bigl|2-{\overline \zeta}({\zeta}+x)\bigr|
\geq 2-|\overline \zeta|(|{\zeta}|+x)\geq  1-x.
\end{equation}
Finally, for $0<\e <1$, $9/10\leq x <1$, under notation \eqref{df:zetaa},
we obtain
\begin{multline*}
|1-{\overline \zeta}x|^2
=(1-tx)^2+4tx\sin^2\frac{\theta}{2}
\leq
\left( 1-\bigl(x-\e(1-x)\bigr) x\right)^2
+4tx\Bigl( \frac{\e(1-x)}{x}\Bigr)^2	\\
\leq(1-x)^2\bigl((1+x+\e x)^2+4t\e^2/x\bigr)\leq
(1-x)^2(9+4/x)\leq 15(1-x)^2.
\end{multline*}
These tree estimates along with \eqref{df:nrep} give  right away \eqref{est:b+}.

The lower estimate \eqref{est:fr2} is trivial.
\end{proof}
\begin{proposition}\label{pr:intx} If\/ $x\notin \Delta_{\e}(x)$, then  
\begin{equation}\label{es:intx}
\bigl(\overline b_2(\zeta , 0)-\overline b_2(\zeta , x)\bigr)^+	
\leq \frac{12}{\e} \frac{(1-|\zeta|)^2}{(1-x)^2} \, .
\end{equation}
\end{proposition}
\begin{proof} It is follows from the representation \eqref{df:HBa}
that
\begin{multline*}
\bigl(\overline b_2(\zeta , 0)-\overline b_2(\zeta , x)\bigr)^+	=
\log^+\left|\dfrac{1-(1-|\zeta|^2)^2}{1-\left(\dfrac{1-|\zeta|^2}{1-\overline{\zeta}x}\right)^2}\right|\\
=\log^+\left|1+\frac{\left(\dfrac{1-|\zeta|^2}{1-\overline{\zeta}x}\right)^2-(1-|\zeta|^2)^2}{1-\left(\dfrac{1-|\zeta|^2}{1-\overline{\zeta}x}\right)^2}\right|
\leq 
\left|\frac{\left(\dfrac{1-|\zeta|^2}{1-\overline{\zeta}x}\right)^2-(1-|\zeta|^2)^2}{1-\left(\dfrac{1-|\zeta|^2}{1-\overline{\zeta}x}\right)^2}\right|\\
=(1-|\zeta|^2)^2 \, \frac{\bigl|1-(1-\overline{\zeta}x)^2\bigr|}{\bigl|(1-\overline{\zeta}x)^2-(1-|\zeta|^2)^2\bigr|}
=(1-|\zeta|)^2 \,\frac{(1+|\zeta|)^2x|2-\overline{\zeta}x|}{|\zeta-x|
\bigl|2-|\zeta|^2-\overline{\zeta}x\bigr|} \, .
\end{multline*}
Hence, if we use \eqref{est:2zet}, and the condition $|\zeta -x|\geq \e(1-x)$, then we obtain
\begin{equation*}
\bigl(\overline b_2(\zeta , 0)-\overline b_2(\zeta , x)\bigr)^+	\leq
(1-|\zeta|)^2 \,\frac{2^2\cdot 3}{\e(1-x)
(1-x)} \, , 	
\end{equation*}
and the estimate \eqref{es:intx} is proved.
\end{proof}
\subparagraph{The upper estimate of \eqref{form:Qb2x}.}
It is follows from the inclusion \eqref{est:D2x} of Proposition
\ref{prop:D2x} that
\begin{multline*} 
Q(x)\leq \int_{ {\Box}_6(x)} \bigl(
\overline b_2(\zeta , 0)-\overline b_2(\zeta , x)\bigr)^+  \d {\nu}_M (\zeta ) \\ \overset{\eqref{df:delex}}{=}
\left( \int_{\Delta_{\e}(x)}+\int_{{\Box}_6(x)\setminus \Delta_{\e}(x)}\right)\bigl(
\overline b_2(\zeta , 0)-\overline b_2(\zeta , x)\bigr)^+  \d {\nu}_M (\zeta )=:I_{\e}(x)+J(x) 
\end{multline*}
for all $9/10\leq x< 1$, where the integrals $I_{\e}(x)$   and $J(x)$ are defined by the last equality. 
 For $I_{\e}(x)$, by estimate \eqref{est:b+} of Proposition \ref{pr:dexes}, we obtain the inequality
\begin{multline*}
I_{\e}(x):=\int_{\Delta_{\e}(x)}\bigl(
\overline b_2(\zeta , 0)-\overline b_2(\zeta , x)\bigr)^+  \d {\nu}_M (\zeta ) \\	\leq 
\int_{\Delta_{\e}(x)} g_{\Delta_{\e}(x)}(\zeta , x)\d {\nu}_M (\zeta ) +\int_{\Delta_{\e}(x)}\log \frac{30}{\e} \d {\nu}_M (\zeta )
\\ \overset{\eqref{est:fr2}}{\leq}
\int_{\Delta_{\e}(x)} g_{\Delta_{\e}(x)}(\zeta , x)\d {\nu}_M (\zeta )
+\frac1{1-\e} \log \frac{30}{\e} 
\cdot \int_{\Delta_{\e}(x)}\frac{(1-|\zeta|)^2}{(1-x)^2}\d {\nu}_M (\zeta )
\end{multline*}
for all $9/10\leq x< 1$. But by the Poisson--Jensen formula
for the disk  $\Delta_{\e}(x)$ we have
\begin{equation*}
\int_{\Delta_{\e}(x)} g_{\Delta_{\e}(x)}(\zeta , x)\d {\nu}_M (\zeta )=
A_{M, \e}(x)-M(x).	
\end{equation*}
Therefore,
\begin{equation}\label{est:intro1}
I_{\e}(x)\leq A_{M, \e}(x)-M(x)+C_{\e}' 
\cdot	
\frac{1}{(1-x)^2}\int_{\Delta_{\e}(x)}{(1-|\zeta|)^2}\d {\nu}_M (\zeta )
\end{equation}
for all $9/10\leq x< 1$ where  the constant $\frac1{1-\e} \log \frac{30}{\e}=:C_e'$ depends only on $\e$. 
 
For the integral $J(x)$, by estimate \eqref{es:intx} of Proposition \ref{pr:intx}, we have
\begin{multline}\label{est:intro}
J(x):=\int_{{\Box}_6(x)\setminus \Delta_{\e}(x)}\bigl(
\overline b_2(\zeta , 0)-\overline b_2(\zeta , x)\bigr)^+  \d {\nu}_M (\zeta )\\
\leq \frac{12}{\e} \int_{{\Box}_6(x)\setminus \Delta_{\e}(x)}\frac{(1-|\zeta|)^2}{(1-x)^2}	\d {\nu}_M (\zeta ).
\\=C_{\e}''\cdot \frac{1}{(1-x)^2}\int_{{\Box}_6(x)\setminus \Delta_{\e}(x)}{(1-|\zeta|)^2}	\d {\nu}_M (\zeta )
\end{multline}
where the constant $12/\e =:C_{\e}''$ also depends only on $\e$.
The addition of \eqref{est:intro1} and \eqref{est:intro}
gives
\begin{equation}
I_{\e}(x)+J(x)\overset{}{\leq }A_{M, \e}(x)-M(x)+C_{\e} 
b_M^{[6]}(x)
\end{equation}
for all $9/10 \leq x <1$ where the constant $C_{\e}:=\max\{C_{\e}', C_{\e}''\}$ depends only on $\e$. Thus Lemma \ref{lem:Qest}
is proved with the constants $a=9/10$ and $\alpha'=6$.$\,\bullet$
\begin{remark} Since $\alpha=\alpha'+1$ in subsection \ref{subsec:alpha},
it follows that we can put $\alpha =7$ in the definition of space \eqref{sp:defq} in Theorem \ref{th:MB}. 
\end{remark}
\section{The proof of Theorem \ref{th:o}}
\setcounter{equation}{0}
We prove Theorem \ref{th:o} as a corollary of Theorem \ref{th:MB}.

By condition of  Theorem \ref{th:o} the function $M\bigm|_{[0, 1)}$ is the increasing positive convex function of\/ $\log$ on $(0, 1)$. Therefore
the function $M(z)=M(|z|)$ is continuous subharmonic function on $\D$, and there exists the positive left-hand derivative $M_-'$  of $M(t)$, $t\in (0,1)$. 
Besides, the function $tM_-'(t)$ is increasing on $(0,1)$.
In particular,
\begin{subequations}\label{equ:M'}
\begin{align}
\int_{r}^{1-}(1-t) \d M(t)&=\int_{r}^{1-}(1-t)  M_-'(t) \d t 
\tag{\ref{equ:M'}i}\label{equ:M'i}\\
\geq 
rM_-'(r) \int_r^{1-}\frac{1-t}{t} \d t &\geq rM_-'(r)\, \frac1{r}\, \frac12 \,(1-r)^2=\frac12 \, M_-'(r)(1-r)^2
\tag{\ref{equ:M'}r}\label{equ:M'r}
\end{align}
\end{subequations}
whence, according to the condition \eqref{cond:inM}, we obtain
\begin{equation}\label{lim:M'}
\lim_{r\to 1-}M_-'(r)(1-r)^2=0. 	
\end{equation}
An easy calculation of Laplacian of $M$ gives the expression  
\begin{equation}\label{df:mesMr}
\d \nu_M(z)=\frac1{2\pi}\d \theta \otimes \d (tM_-'(t)),
\quad z :=te^{i\theta}, \; 0\leq t <1,	
\end{equation}
for the density of Riesz measure $\nu_M$ of $M \in SH(\D)$ (in the sense of distribution theory). 

Suppose that $\alpha$ is the absolute constant for the space \eqref{sp:defq} in Theorem \ref{th:MB}. In our case we have ($z=re^{i\theta}$)
\begin{multline*}
b_M^{[\alpha]}(z)\overset{\eqref{df:bMe}}{=} 
b_M^{[\alpha]}(r)\\
\overset{\eqref{df:box}, \, \eqref{df:mesMr} }{=}
\frac{1}{(1-r)^2} \, 
\frac1{2\pi}\int_{-\alpha(1-r)}^{\alpha (1-r)}\left(\int_{r-\alpha (1-r)}^{1-} (1-t)^2 \d (tM_-'(t))\right) \d \theta \\
=\frac{\alpha}{\pi}\, \frac{1}{1-r}
\int_{r-\alpha (1-r)}^{1-} (1-t)^2 \d (tM_-'(t)).
\end{multline*}
Hence integration by parts gives (taking into account \eqref{lim:M'}
and the positivity of the function $tM_-'(t)$)
\begin{multline*}
b_M^{[\alpha]}(z)\leq 
\frac{\alpha}{\pi}\, \frac{1}{1-r}
\; 2\int_{r-\alpha (1-r)}^{1-}t	M_-'(t)(1-t) \d t\\
\frac{2\alpha}{\pi}\, \frac{1}{1-r}
\, \int_{r-\alpha (1-r)}^{1-}	M_-'(t)(1-t) \d t
\overset{\eqref{equ:M'i}}{=}
\frac{2\alpha}{\pi}\, \frac{1}{1-r}
\, \int_{r-\alpha (1-r)}^{1-}	(1-t) \d M(t)
\\
\overset{\eqref{cond:inM}}{\leq}
\frac{2\alpha}{\pi}\, \frac{1}{1-r} \cdot C\bigl(1-(r-\alpha (1-r))\bigr)
=\frac{2\alpha}{\pi}\, C_M(1+\alpha ),
\end{multline*}
where a constant $C_M$ is dictated by condition \eqref{cond:inM}.
In particular, the condition \eqref{cond:mea} is fulfilled. 
Thus, under the conditions of Theorem \ref{th:o} the space
$\Hol \bigl(\D ; 2A_{M, \e}-M+C_{\e} \cdot b_M^{[\alpha]}\bigr)$	
from \eqref{sp:defq} {\it coincides with the space\/} $\Hol (\D , 2A_{M, \e}-M)$
for every $0< \e <1$.

Let us fix a number $\e \in (0,1)$. The function $M$ is increasing and radial. Therefore,
\begin{equation}\label{es:leA}
A_{M,\e}(z)\leq M\bigl(r+\e (1-r)\bigr), \quad r=|z|.	
\end{equation}
By the mean value theorem there exists a point $r'$
such that  $r\leq r'\leq r+\e (1-r)$ and
\begin{multline*}
M\bigl(r+\e (1-r)\bigr)-M(r)\leq M_-'(r')\cdot \e (1-r)=
M_-'(r')(1-r')\cdot \e \frac{1-r}{1-r'}\\
\leq \frac{\e}{1-\e} \, M_-'(r')(1-r')
\overset{\eqref{equ:M'}}{\leq} \frac{\e}{1-\e} 
\cdot \frac{2}{1-r'}\int_{r'}^{1-}(1-t) \d M(t)
\overset{\eqref{cond:inM}}{\leq} 
\frac{2\e}{1-\e} \cdot C_M,	
\end{multline*}
where $C_M$ is a constant which exists according to the condition \eqref{cond:inM}. Hence, in view of \eqref{es:leA}, we
obtain
\begin{equation*}
2A_{M,\e}(z)-M(z)\leq 2M(r)+2\frac{2\e}{1-\e} \cdot C_M -M(r)=
M(r)+C, \quad r=|z|,	
\end{equation*}
where $C=4C_M\e/(1-\e)$ is a constant. It means that the space 
$\Hol (\D , 2A_{M, \e}-M)$ {\it coincides with the space\/} $\Hol (\D , M)$.
Hence Theorem \ref{th:o} is a special case of Theorem \ref{th:MB}. $\, \bullet$

\section{The case of the uniform Bergman spaces}\label{sec:6} 
\setcounter{equation}{0}

For $V\in {\mc P_{AS}}(\D)$, we put
\begin{equation}\label{char:BC}
\widehat\kappa (V):=\int_0^{1-} \left(\frac1{2\pi}\int_0^{2\pi}
V(te^{i\theta}) \d \theta \right)\, \frac{\d t}{(1-t)^2}.	
\end{equation}
  
\subsection{Weak analogs of Korenblum--Seip's conditions}\label{subsec:KS}
\begin{theorem}\label{th:i} Let $\Lambda=\{ \lambda_k\}$, $k=1, 2, \dots$, $0\notin \Lambda$, be a sequence of points in $\D$, 
and\/ $0\leq p < +\infty$. The following four statements are equivalent:
\begin{enumerate}
\item[{\rm [$\Lambda$]}] the sequence $\Lambda$ is a zero set for the space $A^{-p}$;
\item[{\rm [G]}] there exists a constant $a<1$ such that
$$\sup\limits_{D\in \mc U_0^d(\D ; D(a))}\left(\sum\limits_k    
g_D(\lambda_k, 0) -p \, \widehat\kappa \bigl(g_D(\cdot , 0)\bigr) 
\right) <\infty \, ;$$
\item[{\rm [J]}]
$\sup\limits_{V\in {\mc P_J}(\D)}\left(\sum\limits_k  
V(\lambda_k) -p \, \widehat\kappa (V)\right) <\infty$;
\item[{\rm [AS]}]
$\sup\limits_{V\in {\mc P_{AS}}(\D)}\left(\sum\limits_k   
V(\lambda_k) -p \, \widehat\kappa (V)\right) <\infty$.
\end{enumerate}
\end{theorem}

Evidently, the equivalence [$\Lambda$]$\Longleftrightarrow$[G] coincides with the equivalence $\eqref{i0} \Longleftrightarrow\eqref{ii0}$ of Theorem \ref{th:0}.

\begin{proof} Here we put
\begin{equation}\label{Ma}
M_{p}(z):=p \log \frac{1}{1-|z|}\,	, \quad z\in \D .
\end{equation}
This function is radial and positive. Besides, the restriction 
$M\bigm|_{[0, 1)}$ is the increasing continuous convex function of\/ 
$\log$ on $(0, 1)$, and 
\begin{equation}\label{calk:m'}
\d M_p(t)=\frac1{1-t} \d t,\quad \d\bigl(tM_p'(t)\bigr)=\frac{1}{(1-t)^2} \d t.
\end{equation} 
In particular, condition \eqref{cond:inM} holds,  and $M_p$ is subharmonic on $\D$ with the Riesz measure
\begin{equation}\label{na}
\d \nu_{M_p}	(z)\overset{\eqref{df:mesMr}}{=} \frac{p }{2\pi}\d \theta \otimes \frac{\d t}{(1-t)^2} \, , \quad z=te^{i\theta}\in \D ,\; 0\leq r<1.
\end{equation}
By definition \eqref{sp:HM} we have   
$\Hol(\D , M_{p})=A^{-p}$. In this case condition [G]
coincides with condition \eqref{in:maz1} of Theorem \ref{th:o} with a constant $C$ independent	of $D$. By item ($\mathrm Z_{r}$) of Theorem  \ref{th:o} condition [G] implies [$\Lambda$] and vice versa. 

The evident inclusions 
of the class of all extended Green functions $g_D(\cdot , 0)$, $D\Subset \D$, into $\mc P_{J}(\D) \subset \mc P_{AS} (\D)$ 
give the implications [AS]$\Longrightarrow$[J]$\Longrightarrow$[G].

For the proof of implication [$\Lambda$]$\Longrightarrow$[AS] we use Proposition \ref{pr:dirmajs}.

Let $f_{\Lambda} \in A^{-p}$ be a function with zero set $\Lambda$. It means that the function $\equiv 0$ is harmonic minorant for the  difference
$M_p-\log |f_{\Lambda}|$. Therefore, by Proposition \ref{pr:dirmajs}, the inequalities \eqref{c:ihahM} hold for all for all Arens--Singer functions $V $ on $\D$ with $u=\log |f_{\Lambda}|$,
$\nu_u=n_{\Lambda}$, and $\nu_M=\nu_{M_p}$  In view of 
\eqref{na}, the system of inequalities \eqref{c:ihahM} is equivalent to 
[AS].
  \end{proof}

\subsection{Variants of Luecking's condition}\label{subsec:Lu}

\begin{theorem}\label{th:l} Let $\Lambda=\{ \lambda_k\}$, $k=1, 2, \dots$, $0\notin \Lambda$, be a sequence of points in $\D$ 
and\/ $0\leq p < +\infty$. The following four 
statements are equivalent:

\begin{enumerate}
	\item[{\rm [$\Lambda$]}] The sequence $\Lambda$ is a zero set for the space $A^{-p}$;
\item[{\rm [GL]}] there exists a constant $a<1$ such that 
$$
\hspace{-10mm}\sup\limits_{D\in \mc U_0^d(\D ; D(a))}\left(
 \sum\limits_k   \bigl( 1-|\lambda_k|^2\bigr)^2
\left(\dfrac{1}{\pi}\int_{\D}  \dfrac{g_D(\zeta , 0)\d m(\zeta )}{|1-{\lambda_k}\overline{\zeta}|^4}\right)
 -p\, \widehat\kappa \bigl(g_D(\cdot ,0)\bigr)
\right) <\infty \, ;$$
\item[{\rm [JL]}] $\sup\limits_{V\in {\mc P_J}(\D)}\left( \sum\limits_k   
\bigl( 1-|\lambda_k|^2\bigr)^2
\Bigl(\dfrac1{\pi}\int_{\D}  \dfrac{V(\zeta)\d m(\zeta)}{|1-{\lambda_k}\overline{\zeta}|^4}\Bigr) -p \, \widehat\kappa (V)\right) <\infty $;
\item[{\rm [ASL]}] $\sup\limits_{V\in {\mc P_{AS}}(\D)}\left( \sum\limits_k  
\bigl( 1-|\lambda_k|^2\bigr)^2
\Bigl(\dfrac1{\pi}\int_{\D}  \dfrac{V(\zeta)\d m(\zeta)}{|1-{\lambda_k}\overline{\zeta}|^4}\Bigr) -p \, \widehat\kappa (V)\right) <\infty$.
\end{enumerate}
\end{theorem}

Evidently, the equivalence [$\Lambda$]$\Longleftrightarrow$[GL] coincides with the equivalence $\eqref{i0} \Longleftrightarrow\eqref{ivl0}$ of Theorem \ref{th:0}.

\begin{proof}[Scheme of the proof] 
The evident inclusions 
of the class of all extended Green functions $g_D(\cdot , 0)$, $D\Subset \D$, into $\mc P_{J}(\D) \subset \mc P_{AS} (\D)$ 
give the implications [ASL]$\Longrightarrow$[JL]$\Longrightarrow$[GL].

Put
\begin{equation}\label{df:kL}
K_{\Lambda}(z)=\frac{|z|^2}{2}\sum_{\lambda_k\in \Lambda}
\frac{\bigl(1-|\lambda_k|^2\bigr)^2}{|1-\lambda_k \overline{z}|^2}.
\end{equation}
 The function $K_{\Lambda}$ is a subharmonic continuous function
on $\D$.
 Easy calculations give
\begin{multline*}
\Delta K_{\Lambda}(z)=4\frac{\partial^2 }{\partial z \partial \overline{z}}
\, K_{\Lambda}(z)\\
=2\sum_{k} \bigl(1-|\lambda_k|^2\bigr)^2
\frac{\partial^2 }{\partial z \partial \overline{z}}
\, \frac{z \overline{z}}{(1-\lambda_k \overline{z})(1-\overline{\lambda}_k z)}
\\=2 \sum_{k} \bigl(1-|\lambda_k|^2\bigr)^2
\frac{1}{|1-{\lambda_k}\overline{z}|^4}\, .
\end{multline*}
Hence we get 
\begin{equation}\label{mes:KL}
\d \nu_{K_{\Lambda}}(z)= \frac1{\pi} \sum_{k} 
\frac{\bigl(1-|\lambda_k|^2\bigr)^2}{|1-{\lambda_k}\overline{z}|^4}  \d m(z).
\end{equation}
By Luecking's criterion \cite[Theorem A, (b)]{Luecking} {\it a sequence $\Lambda=\{ \lambda_k\}$ is zero set for\/
$A^{-p}$ if and only if for $M_{p}$ from \eqref{Ma}, the function
$M_{p}-K_{\Lambda}$ admits a harmonic minorant}.

Let $M=M_{p}$, $\nu_M=\nu_{M_p}$ (see \eqref{na}), and $u=K_{\Lambda}$. It follows from Proposition \ref{pr:dirmajs} the same way as in the proof of implication [$\Lambda$]$\Longrightarrow$[AS] of previous Theorem that 
the implication [$\Lambda$]$\Longrightarrow$[ASL] holds.

For the proof of implication [GL]$\Longrightarrow$[$\Lambda$]
we can use subharmonic analogs of Corollary \ref{th:Mc}
(item [{\rm ($\mathrm Z_{\sic}$)}]), Theorem 
\ref{th:MB} (item [{\rm ($\mathrm Z_{\D}$)}]), and Theorem \ref{th:o}
(item [{\rm ($\mathrm Z_{r}$)}]) with the measure $\nu_{K_{\Lambda}}$ in place of measure $n_{\Lambda}$. In particular, the last is the statement
\begin{enumerate}
	\item[{\rm [$\mathrm S\mathrm H_r$]}] {\it Let the function\/ $M$ be the same as in Theorem\/ {\rm \ref{th:o}}, and\/ $K$ be a subharmonic continuous function on\/ $\D$ with Riesz measure\/ $\nu_K$. The function $M-K$ admits a harmonic minorant if and only if there are constants $a<1$ and $C$   
such that the inequality  
\begin{equation*}
\int_{\D} g_D(\zeta , 0) \d \nu_K(\zeta)\leq
\int_0^{1^-} \left(\frac{1}{2\pi} \int_0^{2\pi} g_D (te^{\theta} , 0 )\d \theta \right)\d \bigl( tM_-'(t)\bigr) +C
\end{equation*}
holds for each domain $D \in \mc U_0^d \bigl(\D ; D(a ) \bigr)$}. 	
\end{enumerate}
The system of inequalities in the statement [$\mathrm S\mathrm H_r$] is equivalent to the statement [GL] with $K_{\Lambda}$ and $M_p$
from \eqref{df:kL} and  \eqref{Ma} in place of $K$ and $M$ respectively according to \eqref{mes:KL} and \eqref{calk:m'}.
It proves the  equivalence [GL]$\Longleftrightarrow$[$\Lambda$].
\end{proof}
\begin{remarks} 1. The summands
in first sums of [GL], [JL], and [ASL] are values at points $\lambda_k$ of the Berezin transforms  (see \cite[Ch.~2]{HKZ}) of Green, Jensen and Arens--Singer functions respectively.

2. In Theorems \ref{th:i} and \ref{th:l} in items [G] and [GL] we can suppose that the domain $D\subset \D$ has final number of points of contact with circumference $\partial \D$ since the integral     
\eqref{char:BC} with the integrand $V=g_D(\cdot , 0)$ 
is a finite number.
\end{remarks}

More detailed statement of this work will be submit in the journal ``Matematicheskii sbornik'' in Russian (the translation from Russian into English is ``Sbornik: Mathematics'').

\end{document}